\documentclass[10pt]{amsart} 

\usepackage{microtype}
\usepackage{amsmath, amssymb, mathrsfs} 
\usepackage{mathtools}
\usepackage[mathscr]{euscript} 
\usepackage{bbm}
\usepackage[breaklinks=true]{hyperref} 
\usepackage{comment}
\usepackage{url}
\usepackage{tikz-cd}
\usetikzlibrary{graphs,decorations.pathmorphing,decorations.markings}

\newlength{\mylength}
\usepackage{enumitem}
\setlist{listparindent=\parindent, itemsep=0cm, parsep=\mylength, topsep=0cm}

\newenvironment{cd}{\begin{equation*}\begin{tikzcd}}{\end{tikzcd}\end{equation*}\ignorespacesafterend}
\newcommand{\e}[1]{\begin{align*} #1 \end{align*}}

\usepackage{amsthm}

\newtheoremstyle{mythm}
{\mylength}
{0pt}
{\itshape}
{0pt}
{}
{. }
{0em}
{\thmnumber{#2. }\bfseries{\thmname{#1}\thmnote{ (#3)}}}

\newtheoremstyle{mydef}
{\mylength}
{0pt}
{}
{0pt}
{}
{. }
{0em}
{\thmnumber{#2. }\bfseries{\thmname{#1}\thmnote{ (#3)}}}

\newtheoremstyle{myrmk}
{\mylength}
{0pt}
{}
{0pt}
{}
{.\ }
{ }
{\thmnumber{#2. }\emph{\thmname{#1}\thmnote{ (#3)}}}

\makeatletter
\renewenvironment{proof}[1][\proofname]{\par
	\pushQED{\qed}%
	\normalfont \topsep6\p@\@plus6\p@\relax
	\noindent\emph{#1.} 
	\ignorespaces
}{%
\popQED\endtrivlist\@endpefalse
}
\makeatother

\theoremstyle{mythm} 
\newtheorem*{thm*}{Theorem} 
\newtheorem*{lem*}{Lemma} 
\newtheorem*{cor*}{Corollary} 
\newtheorem*{claim*}{Claim} 
\newtheorem*{prop*}{Proposition}

\theoremstyle{mydef} 
\newtheorem*{conj*}{Conjecture}
\newtheorem*{que*}{Question}

\theoremstyle{myrmk}
\newtheorem*{rmk*}{Remark} 
\newtheorem*{rmks*}{Remarks} 
\newtheorem*{ex*}{Example}

\makeatletter
\def\theoremhead{\@startsection{subsubsection}{3}%
	\z@{.5\linespacing\@plus.7\linespacing}{-.5em}%
	{\normalfont\bfseries}}
\makeatother

\newenvironment{customtheorem}[1]%
{\theoremhead{#1} \begingroup \itshape}%
{\endgroup \par} 
{\theoremhead{#1} }%
{\par} 
{\subsubsection{Remark}}%
{\par} 

\newenvironment{cor}{\begin{customtheorem}{Corollary}}{\end{customtheorem}}
\newenvironment{thm}{\begin{customtheorem}{Theorem}}{\end{customtheorem}}
\newenvironment{lem}{\begin{customtheorem}{Lemma}}{\end{customtheorem}}
\newenvironment{prop}{\begin{customtheorem}{Proposition}}{\end{customtheorem}}

\newcommand{\on}{\operatorname}
\renewcommand{\setminus}{\smallsetminus} 
\newcommand{\ol}{\overline} 
\newcommand{\ul}{\underline} 

\newcommand{\BZ}{\mathbb{Z}}

\newcommand{\pr}{\on{\ms{pr}}}
\newcommand{\id}{\on{\ms{id}}}

\newcommand{\la}{\langle}
\newcommand{\ra}{\rangle} 
\newcommand{\wt}{\widetilde} 
\newcommand{\wh}{\widehat}
\newcommand{\mb}{\mathbf}
\newcommand{\mf}{\mathfrak}
\newcommand{\mc}{\mathscr}
\newcommand{\mr}{\mathrm}
\newcommand{\ms}{\mathsf} 

\newcommand{\hra}{\hookrightarrow}
\newcommand{\sra}{\twoheadrightarrow}

\newcommand{\Hom}{\on{Hom}}
\newcommand{\Fl}{{\mc{F}\hspace{-.3mm}\ell} }
\newcommand{\curHom}{\mathcal{H}\kern -.5pt om}
\newcommand{\curExt}{\mathcal{E}\kern -.5pt xt}
\newcommand{\EExt}{\on{\mathbb{E}\kern -.5pt \mathrm{xt}}}
\newcommand{\curSpec}{\mathcal{S}\kern -0.5pt pec}
\newcommand{\curProj}{\mathcal{P}\kern -0.5pt roj}

\newcommand{\oh}{\mc{O}}

\DeclareMathOperator*\colim{colim}
\newcommand{\LKE}{\on{\mr{LKE}}}
\newcommand{\Res}{\on{\mr{Res}}}

\newcommand{\Domc}{\ms{Dom}^{\ms{c}}}

\newcommand{\Embc}{\ms{Emb}^{\ms{c}}}
\newcommand{\Domd}{\ms{Dom}^{\ms{d}}}
\newcommand{\Embd}{\ms{Emb}^{\ms{d}}}
\newcommand{\Ic}{\ms{Rig}^{\ms{c}}}
\newcommand{\Id}{\ms{Rig}^{\ms{d}}}

\newcommand{\pt}{\ms{pt}}

\newcommand{\bb}[1]{[\![#1]\!]}

\renewcommand{\mod}{\on{\!-mod}}

\newcommand{\sset}{\mathsf{SSet}}
\newcommand{\nec}{\mathsf{Nec}}

\usepackage{bm}

\newcommand{\op}{{\mr{op}}}

\newcommand{\assoc}{\mathrm{N}(\Delta^{\op})}
\newcommand{\ner}{\mathrm{N}}
\newcommand{\alg}{\mathsf{Alg}}
\newcommand{\cat}{\mathsf{DGCat}}
\renewcommand{\mod}{\text{-}\mathsf{mod}}

\newcommand{\subsubsectiona}{\subsubsection{}\hspace{-0.7em}}

\newcommand{\LG}{\mc{L}G}
\newcommand{\xra}{\xrightarrow} 

\newcommand{\prl}{\mc{P}\mr{r}^{\mr{L}}}

\newcommand\vphi\varphi
\usepackage[margin=1.5in]{geometry}

\title{The affine Hecke category is a monoidal colimit}
\author{James Tao} 
\address{Massachusetts Institute of Technology, Cambridge, MA 02139, USA}
\email{james.h.tao@gmail.com}
\author{Roman Travkin}
\address{Skolkovo Institute of Science and Technology, Moscow, Russia}
\email{roman.travkin2012@gmail.com}
\date{August 21, 2026}

	\makeatletter
	\newcommand{\needspace}[1]{%
		\begingroup
		\setlength{\dimen@}{#1}%
		\vskip\z@\@plus\dimen@
		\penalty -100\vskip\z@\@plus -\dimen@
		\vskip\dimen@
		\penalty 9999%
		\vskip -\dimen@
		\vskip\z@skip 
		\endgroup
	}
	
	\newcommand{\Needspace}{\@ifstar{\@sneedsp@}{\@needsp@}}
	
	\newcommand{\@sneedsp@}[1]{\par \penalty-100\begingroup
		\setlength{\dimen@}{#1}%
		\dimen@ii\pagegoal \advance\dimen@ii-\pagetotal
		\ifdim \dimen@>\dimen@ii
		\break
		\fi\endgroup}
	
	\newcommand{\@needsp@}[1]{\par \penalty-100\begingroup
		\setlength{\dimen@}{#1}%
		\dimen@ii\pagegoal \advance\dimen@ii-\pagetotal
		\ifdim \dimen@>\dimen@ii
		\ifdim \dimen@ii>\z@
		\vfil
		\fi
		\break
		\fi\endgroup}
	\makeatother
	
	\newcommand{\mysubsec}[1]{\needspace{0.8in} \subsection{#1} }
	
\newcommand{\spine}{\ms{spine}}
\newcommand{\Shv}{\mc{D}}
\newcommand{\D}{\mc{D}}

\begin{document}

	\begin{abstract}
	\vspace{-4mm}
		Let $G$ be a semisimple simply-connected algebraic group over an algebraically closed field of characteristic zero. We prove that the affine Hecke category associated to the loop group of $G$ is equivalent to the colimit, evaluated in the $\infty$-category of monoidal stable $\infty$-categories, of the finite type Hecke subcategories associated to standard parahoric subgroups. The main ingredient is an inductive characterization of colimits indexed by (sufficiently nice) bistratified categories. Our method is very general and can be used to prove a number of analogous `colimit theorems,' e.g.\ for $\mc{D}$-modules on the loop group. 
	\end{abstract}

	\maketitle
	
	\vspace{-10mm}
	\setlength{\parskip}{0.8mm}	

	\tableofcontents 
	
	\newpage
	
	\setlength{\parskip}{\mylength}
	
	\section{Introduction} \label{intro} 
		
	\mysubsec{Monoidal colimits} 
	Let $G$ be a simply-connected 	Kac--Moody group with Weyl group $W$, and let $I$ be a set of simple reflections. We say that a subset $J \subset I$ is \emph{finite type} if the subgroup $W_J \subset W$ generated by the simple reflections in $J$ is finite. In this case, we have the standard type-$J$ parabolic subgroup $P_J \subset G$. 
	
	It is an old idea that the following should be true: 
	\[
		G = \colim_{\substack{J \subset I \\ J \text{ fin.\ type}}} P_J, 
	\]
	where the colimit is taken in groups or monoids. This idea is implicit in Tits's original presentation of Kac--Moody groups via generators and relations, see~\cite{tits}. One can show that the statement is equivalent to the simple connectivity of the Tits building of $G$. Working with topological groups gives a stronger statement, which was deduced in~\cite[Thm.\ 4.2.3]{nitu} from the contractibility of a topological analogue of the Tits building. In another paper~\cite{hyper}, we use the contractibility of the Tits building to show that the statement remains true if $G$ and $P_J$ are replaced by the $\infty$-categories of $\D$-modules $\Shv(G)$ and $\Shv(P_J)$, and the colimit is taken in monoidal stable $\infty$-categories (or DG-categories). 
	
	With this in mind, it is natural to conjecture that the statement remains true if $G$ and $P_J$ are replaced by their Hecke categories: 
	\[
		\Shv(B \backslash G / B) = \colim_{\substack{J \subset I \\ J \text{ fin.\ type}}} \Shv(B \backslash P_J / B). 
	\]
	This turns out to be true, and it is the main theorem of the present paper (Theorem~\ref{mainthm}). However, the Tits building is no longer useful, essentially because the left $B$-quotient smears together the points of the flag variety $G/B$ which ought to correspond to different chambers in the Tits building. 
	
	Instead, we take a more direct and combinatorial approach, which essentially handles one Schubert cell at a time. For simplicity, let us return to the statement $G = \colim_J P_J$ and work in the ambient category of schemes. The first step is to rewrite the monoidal colimit as an `amalgamated product': 
	\begin{itemize}
		\item Let $\ms{Words}$ be the category of sequences of finite type subsets $(J_1, \ldots, J_n)$, where a morphism $(J_1, \ldots, J_n) \to (J'_1, \ldots, J'_{n'})$ is a weakly-increasing map $\psi : [1, n] \to [1, n']$ such that, for each $j \in [1, n']$, we have $J'_j \supseteq \bigcup_{i \in \psi^{-1}(j)} J_i$. 
		\item Define a functor $E : \ms{Words} \to \ms{Schemes}$ as follows: 
		\begin{itemize}
			\item $E(J_1, \ldots, J_n) := P_{J_1} \overset{B}{\times} \cdots \overset{B}{\times} P_{J_n}$. (Cf.~\cite[\S 1]{mathieu}.) 
			\item The behavior of $E$ on morphisms comes from the embeddings $P_J \hra P_{J'}$ for $J \subset J'$ and the multiplication maps for each $P_J$. 
		\end{itemize}
		\item Then $\colim_J P_J$ identifies with the (non-monoidal) colimit $\colim_{\ms{Words}} E$. 
	\end{itemize}
	Thus, the problem reduces to computing a (non-monoidal) colimit indexed by $\ms{Words}$. 
	
	\begin{rmks*} \needspace{0.5in} \leavevmode
		\begin{enumerate}[label=(\arabic*)]
			\item The method of the present paper is very general. For example, it proves all of the above monoidal colimit statements, not just the statement for Hecke categories. 
			\item The above statements are useful for constructing actions of $G$, $\Shv(G)$, or $\Shv(B \backslash G / B)$. 		\end{enumerate}
	\end{rmks*}
	
	\mysubsec{Reedy-like categories} Recall that a Reedy category is a category $\mc{C}$ equipped with a rank function $r : \mathrm{Ob}(\mc{C}) \to \BZ_{\ge 0}$ and two subcategories of `upward' and `downward' maps, such that there is a functorial factorization of every map into a downward map followed by an upward map. It is well-known that colimits indexed by Reedy categories can be computed inductively, moving from lower to higher ranks. 
	
	\begin{rmk*}
		The main example of a Reedy category is the simplex category $\mb{\Delta}$. The rank function sends each simplex to its dimension: $\Delta^n \mapsto n$. The upward and downward maps are given by injective and surjective maps of finite nonempty ordered sets, respectively. 
	\end{rmk*}
	
	Our main idea is to think of $\ms{Words}$ as a Reedy category, via the following choices: 
	\begin{itemize}
		\item Define the functor $D : \ms{Words} \to \ms{Varieties}$ via $D(J_1, \ldots, J_n) = E(J_1, \ldots, J_n) / B$. 
		
		This is a \emph{Demazure variety}, i.e.\ a twisted product of finite-type flag varieties. 
		\item A map is upward if its $D$-image is a closed embedding of varieties. 
		\item A map is downward if its $D$-image is a dominant map of varieties. 
		\item The functorial factorization in $\ms{Words}$ corresponds to the factorization
		\[
			D(J_1, \ldots, J_n) \to (\text{scheme-theoretic image}) \to D(J'_1, \ldots, J'_{n'}), 
		\]
	\end{itemize}
	Unfortunately, this idea suffers from two problems: 
	\begin{itemize}
		\item The scheme-theoretic image may fail to lie in the image of $D$, in which case we cannot get a factorization of $(J_1, \ldots, J_n) \to (J'_1, \ldots, J'_{n'})$. 
		\item $r$ fails to be strictly decreasing along all nonidentity downward maps. Indeed, $D$ sends some nonidentity downward maps to \emph{birational} maps of varieties. 
	\end{itemize}
	
	The first problem is easy to resolve. Generalize the notion of `Demazure variety' to include twisted products of finite-type Schubert varieties. The aforementioned scheme-theoretic image is always of this form. Passing to this larger category $\ms{Words}'$ does not change the colimit, essentially because every Schubert variety embeds into a flag variety. 
	
	The second problem is more serious. Fortunately, there is already a generalization of the notion of a Reedy category, due to Shulman~\cite{reedy}, which allows for the possibility that some nonidentity downward maps fail to decrease the rank. 
	
	\mysubsec{Outline} We can now summarize the contents of the present paper: 
	\begin{itemize}
		\item In Section~\ref{s-bis}, we enhance a small part of Shulman's theory, so as to make it suitable for computing the colimit of a functor from a Reedy-like category to an $\infty$-category. We also make two observations which are specific to our application: 
		\begin{itemize}
			\item In~\ref{extra-bis}, we observe that computing colimits is especially easy if certain categories of downward maps are contractible. 
			\item In~\ref{rigidcat}, we explain how to present a Reedy-like category by specifying the upward maps, downward maps, and certain commutation relations between them. Most Reedy-like categories cannot be presented in this way, but our framework applies to $\mb{\Delta}$, $\ms{Varieties}$, $\ms{Schemes}$, and $\ms{Words}'$. 
		\end{itemize}
		\item In Section~\ref{dem}, we show that $\ms{Words}'$ is Reedy-like and satisfies the aforementioned condition on the contractibility of downward maps. 
		\item In Section~\ref{s-main}, we combine the previous ideas and prove the main theorem. Nothing nontrivial happens in this section, but we have taken pains to rigorously treat the $\infty$-categorical aspects of the following topics: monoidal structures, monoidal colimits, and gluing (a.k.a.\ recollement) of $\D$-module categories. For the last topic, we use a very convenient framework for gluing $\infty$-categories which is developed in~\cite{strat}.
	\end{itemize}
	In the sequel, $\ms{Words}$ is denoted $(\mc{P}_{I, \mr{fin}})^{\amalg}_{\mr{act}}$, see~\ref{ss-amalg}. Also, $\ms{Words}'$ is denoted $\ms{Rig}^{\ms{d}}$, see~\ref{concrete}. 
	
	\mysubsec{Notations and conventions} 
	Let $[n] := \{0, 1, \ldots, n\}$. 
	
	A functor is called \emph{initial} or \emph{final} if precomposition by it preserves limits or colimits, respectively. (It is perhaps more common to use final/cofinal in place of initial/final.) 
	
	The \emph{comma category} associated to a pair of functors 
	\begin{cd}
		\mc{C}_1 \ar[r, "F_1"] & \mc{D} \ar[r, leftarrow, "F_2"] & \mc{C}_2
	\end{cd}
	is usually denoted $(F_1 \downarrow F_2)$. However, we prefer to use angle brackets: $\la F_1 \downarrow F_2 \ra$. 
	(Its objects are triples $(c_1,c_2,f)$ where $c_i\in \mc C_i$ and $f\colon F_1(c_1)\to F_2(c_2)$, with morphisms given by $\Hom_{\la F_1 \downarrow F_2\ra}((c_1,c_2,f),(c_1',c_2',f')) = \{(\vphi_1,\vphi_2)\in\Hom(c_1,c_1')\times\Hom(c_2,c_2')\mid F_2(\vphi_2)\circ f = f' \circ F_1(\vphi_1)\}$.)
	There are evident forgetful functors $\ms{tail}\colon\la F_1\downarrow F_2\ra \to \mc C_1$ and  $\ms{head}\colon\la F_1\downarrow F_2\ra \to \mc C_2$. When the functors {$F_1,F_2$} are clear from context, it will also be convenient to denote this category as $\la \mc{C}_1 \xra{\mc{D}} \mc{C}_2 \ra$. If $F_1$ is the embedding of a single object $d$, and $F_2$ is the identity functor, the resulting \emph{undercategory} is also denoted $\mc{D}_{d/}$ or $\la d \to \mc{D} \ra$. Similarly, overcategories are denoted $\mc{D}_{/d}$ or $\la \mc{D} \to d \ra$. 
	
	An \emph{$\infty$-category} is an $(\infty, 1)$-category modeled as a weak Kan complex, as in~\cite{htt} and \cite{ha}. For convenience, we do not distinguish between a category and its nerve, so that any category \emph{is} an $\infty$-category. In this paper, all categories are small, but some $\infty$-categories are not small. 
	
	A \emph{categorical equivalence} of simplicial sets is a weak equivalence in the Joyal model structure, while a \emph{homotopy equivalence} is a weak equivalence in the Kan model structure. (See~\cite[Def.\ 1.1.5.14]{htt}.) When working with $\infty$-categories, we abbreviate `categorical equivalence' as `equivalence,' but the phrase `homotopy equivalence' will never be abbreviated. For us, \emph{contractible} means `weakly contractible' (i.e.\ homotopy equivalent to a point). 
	
	$\mathsf{Spaces}$ is the $\infty$-category of topological spaces (or $\infty$-groupoids), and $\ms{SSet}$ is the category of simplicial sets. 
	
	Let $\prl$ be the $\infty$-category whose objects are presentable $\infty$-categories and whose morphisms are functors which preserve small colimits, equipped with Lurie's tensor product~\cite[4.8.1]{ha}. Let $\ms{Vect} \in \alg(\prl)$ be the unbounded derived category of vector spaces over a field of characteristic zero~\cite[1.5.3]{ha}, and let $\cat$ be the $\infty$-category of $\ms{Vect}$-modules in $\prl$ which are stable. (See~\cite[Vol.\ 1, Chap.\ 1]{gr}.) We let $\Shv(-)$ denote the $\infty$-category of $\D$-modules as crystals developed in~\cite[Vol.\ 2, Chap.\ 4]{gr}, so that $\Shv(-)$ is a functor from finite type stacks to $\cat$. Here are two facts about this framework which we will need: 
	\begin{itemize}
		\item \cite[Cor.\ 5.5.3.4]{htt} implies that a colimit in $\cat$ can be computed as the limit over the diagram obtained by passing to right adjoints. 
		\item $\mc{D}$-modules satisfy descent along blow-up squares of finite type algebraic stacks. (See the proof of~Lemma~\ref{blowup}.) 
	\end{itemize}

	\mysubsec{Relation to previous versions}
	The present version supersedes versions 1--3 of this paper (\texttt{arXiv:2009.10998}), which proved the main theorem by a different and more combinatorial method, organized around a homotopical analogue of the Coxeter deletion property for the `convolution Schubert $1$-category' $\ms{Word}^{1}_{\mr{fr}}$. That framework does not appear in the present version, having been replaced by the Demazure category of Section~\ref{dem}. Readers who arrive here by following a reference to the earlier framework, such as~\cite[App.\ C]{LNY}, should consult \texttt{arXiv:2009.10998v3}. 

	\mysubsec{A related work} 	
	After completion of an earlier draft of this paper, we were informed of a project by Li, Nadler, Varshavsky and Yun achieving related results by a different method.
	
	\mysubsec{Applications}
	The affine Hecke category and its cocenter play an important role in geometric Langlands duality and related areas. The colimit presentation established in the present paper yields new effective tools for their study; it is used in a crucial way in~\cite{LNY}, where it serves as a starting point for a geometric description of the cocenter of the affine Hecke category. See also~\cite{zhu-tame}, which uses the companion result of~\cite{hyper}; that result also follows from the method of the present paper.
	
	One of the original motivations for this project was a potential construction of a functor from the affine to the finite type A Hecke category, extending the construction of~\cite{kostya-thesis}.  However, currently we do not know how to construct all the higher homotopies in the relevant diagram of monoidal functors.
	
	\mysubsec{Acknowledgments} 
	We are indebted to Roman Bezrukavnikov, Kostya Tolmachov, and Sam Raskin for helpful conversations which took place during the course of this project. In addition, we would like to thank Dennis Gaitsgory, Penghui Li, and David Nadler for their cordial correspondence regarding the aforementioned work, and Yakov Varshavsky for generously providing a detailed overview of his results on loop group invariants. The first author was supported by the NSF GRFP, grant no.\ 1122374.
	
	\clearpage
	\needspace{4em}
	\section{\texorpdfstring{$\infty$}{oo}-bistratified categories and colimits} \label{s-bis}  
	
	\mysubsec{Joyal fibrant replacement via necklaces} \label{ss-neck} 
	
	\subsubsectiona  \label{neckthm} 
	By definition, a fibrant replacement functor for the Joyal model structure gives a `universal' way to turn simplicial sets into $\infty$-categories. In~\cite{rigid}, Dugger and Spivak defined a very convenient fibrant replacement functor, which we summarize below. 
	
	First, define a \emph{necklace} to be a nonempty simplicial set of the form 
	\[
		\Delta^{n_0} \vee \Delta^{n_1} \vee \cdots \vee \Delta^{n_k}, 
	\]
	where each wedge means that the vertex $n_i \in \Delta^{n_i}$ is glued to the vertex $0 \in \Delta^{n_{i+1}}$. Each necklace has a unique `first vertex' $0 \in \Delta^{n_0}$ and `last vertex' $n_k \in \Delta^{n_k}$. Its \emph{length} is $n_0 + \cdots + n_k$. The point $\Delta^0$ is the unique necklace whose first and last vertices coincide. 
		
	Let $\nec$ be the category whose objects are necklaces and whose morphisms are maps of simplicial sets which send the first vertex to the first vertex and the last vertex to the last vertex. For any simplicial set $S$ and vertices $a, b \in S_0$, we define $\nec_S(a, b)$ to be the category of pairs $(T, \pi)$ where $T \in \nec$ and $\pi : T \to S$ is a map sending the first vertex to $a$ and the last vertex to $b$. As usual, morphisms in $\nec_S(a, b)$ are commutative triangles. 
	
	\begin{thm} \label{neck1} 
		Let $S$ be a simplicial set, let $a, b \in S_0$ be vertices, and let $S \to \wt{S}$ be a fibrant replacement in the Joyal model structure. Then the space of maps $\Hom_{\wt{S}}(a, b)$ is homotopy equivalent to $\nec_S(a, b)$. 
	\end{thm}
	\begin{proof}
		Theorem~1.3 in~\cite{rigid} states that there is a zig-zag of weak equivalences between two simplicial categories $\mf{C}^{\mr{nec}}(S)$ and $\mf{C}[S]$. The simplicial category $\mf{C}^{\mr{nec}}(S)$ is defined so that 
		\[
		\Hom_{\mf{C}^{\mr{nec}}(S)}(a, b) = \on{N}(\nec_S(a, b)),
		\]
		while the simplicial category $\mf{C}[S]$ is defined in~\cite[1.1.5]{htt} and has the property that the topological nerve of $\lvert\mf{C}[S]\rvert$ is a fibrant replacement for $S$. 
	\end{proof}

	\subsubsectiona  \label{neckreform} 
	We are interested in the case when $S$ is a simplicial subset of a category $\mc{C}$. In this case, the definition of $\nec_S(a, b)$ can be reformulated as follows: 
	\begin{itemize}
		\item For any simplex $\Delta^n$, let $\spine(\Delta^n)$ be the smallest simplicial subset which includes the edges $\Delta^{\{i-1, i\}}$ for $i \in [1, n]$. 
		\item An object of $\nec_S(a, b)$ is a triple $(\Delta^n, \phi, P)$, where $\phi : \spine(\Delta^n) \to S$ is any map, and $P \subseteq [n]$ is any subset which contains $0$ and $n$, such that the following holds: 
		\begin{itemize}
			\item For any consecutive elements $p, q \in P$, note that $\spine(\Delta^{[p, q]}) \xra{\phi} S \hra \mc{C}$ extends uniquely to a map $\Delta^{[p, q]} \to \mc{C}$. We require that the latter factors through $S$. 
		\end{itemize}
		\item A morphism $(\Delta^{n_1}, \phi_1, P_1) \to (\Delta^{n_2}, \phi_2, P_2)$ is a map $F : \Delta^{n_1} \to \Delta^{n_2}$ such that 
		\begin{itemize}
			\item $F(P_1) \supseteq P_2$. 
			\item For every $i \in [1, n_1]$, the $\mc{C}$-composite of the arrows $\spine(\Delta^{[F(i-1), F(i)]}) \xra{\phi_2} S$ equals the arrow $\Delta^{\{i-1, i\}} \xra{\phi_1} S$. 
		\end{itemize}
	\end{itemize}
	This is equivalent to the original definition of $\nec_S(a, b)$ because $(\Delta^n, P)$ encodes a length-$n$ necklace, where $P$ specifies the gluing points. A map from this necklace to $S$ is completely determined by its restriction to $\spine(\Delta^n)$, because composition of $\mc{C}$-morphisms is unique.

	\mysubsec{\texorpdfstring{$\infty$}{oo}-bistratified categories} \label{ss-bis} 
	The following is a generalization of the notion of a Reedy category. We learned this material from~\cite[\S 6]{reedy}. Our treatment is the same, except we have replaced the word `connected' with the word `contractible,' and our notion of `categorical equivalence' is defined in the context of $\infty$-categories. 
	
	\subsubsectiona 
	Let $\mc{C}$ be a small category, let $J$ be an ordered set, and fix a map $r : \on{\ms{Obj}}(\mc{C}) \to J$. The following definitions pertain to a morphism $f : c_1 \to c_2$ in $\mc{C}$. 
	\begin{itemize}
		\item If $r(c_1) < r(c_2)$, then $f$ is \emph{increasing}. 
		\item If $r(c_1) > r(c_2)$, then $f$ is \emph{decreasing}. 
		\item If $r(c_1) = r(c_2)$, then $f$ is \emph{level}. 
		\item Let $\ms{Fact}_f$ be the category whose objects are factorizations of $f$ as 
		\begin{cd} 
			c_1 \ar[r] & c \ar[r] & c_2
		\end{cd} 
		where $r(c_1) > r(c) < r(c_2)$, and morphisms are commutative diagrams 
		\begin{cd}
			c_1 \ar[r] \ar[d, dash, shift left = 0.5] \ar[d, dash, shift right = 0.5] & c \ar[r] \ar[d] & c_2 \ar[d, dash, shift left = 0.5] \ar[d, dash, shift right = 0.5] \\
			c_1 \ar[r] & c' \ar[r] & c_2
		\end{cd}
		where the top and bottom rows are factorizations as above. 
		\item If $\ms{Fact}_f$ is empty, we say that $f$ is \emph{basic}. 
	\end{itemize}
	We say that $\mc{C}$ is \emph{$\infty$-bistratified with rank function $r$} if the following are true: 
	\begin{enumerate}[label=(B\arabic*)]
		\item The basic level morphisms form a subcategory of $\mc{C}$, i.e.\ all identity morphisms are basic level, and the composite of two basic level morphisms is basic level. \label{cond-b1}
		\item For any level morphism $f$ which is not basic, $\ms{Fact}_f$ is contractible. \label{cond-b2}
	\end{enumerate}
	In this case, every $j \in J$ determines three subcategories $\mc{C}_{<j}, \mc{C}_{\le j}, \mc{C}_{j} \subset \mc{C}$ as follows. 
	\begin{itemize}
		\item $\mc{C}_{<j} \subset \mc{C}$ is the full subcategory such that $c \in  \mc{C}_{<j}$ if and only if $r(c) < j$. 
		\item $\mc{C}_{\le j} \subset \mc{C}$ is a full subcategory which is defined similarly. 
		\item $\mc{C}_j \subset \mc{C}$ is the non-full subcategory whose objects are $c \in  \mc{C}$ such that $r(c) = j$ and whose morphisms are the basic level morphisms. 
	\end{itemize}
	Also, define the \emph{$j$-th bigluing datum} to be the simplicial subset $\mc{C}_{\le j}^\circ \subset \mc{C}_{\le j}$ which contains all objects, all morphisms except for non-basic morphisms satisfying $r(c_1) = r(c_2) = j$, and all possible simplices of dimension $\ge 2$ (whose edges are all in the above set of morphisms). 
	
	\begin{thm}\label{recognize}
		The embedding $\mc{C}_{\le j}^\circ \hra \mc{C}_{\le j}$ is a categorical equivalence. 
	\end{thm}
	
	This theorem corresponds to~\cite[Thm.\ 6.13]{reedy}. Its proof occupies the rest of this subsection.

	\subsubsection{Necklaces} 
	Factor the embedding $\mc{C}_{\le j}^\circ \hra \mc{C}_{\le j}$ as follows: 
	\begin{cd}[column sep = 1.1in] 
		\mc{C}_{\le j}^\circ \ar[r, hookrightarrow, "\text{fibrant replacement}"] &  (\mc{C}_{\le j}^\circ)_{\infty} \ar[r] & \mc{C}_{\le j}
	\end{cd}
	The first map is a categorical equivalence by construction, so we only need to show that the second map is a categorical equivalence. It suffices to show that, for all $c_1, c_2 \in \mc{C}_{\le j}$, the map 
	\[
	\rho : \Hom_{(\mc{C}_{\le j}^\circ)_{\infty}}(c_1, c_2) \to \Hom_{\mc{C}_{\le j}}(c_1, c_2)
	\]
	is a homotopy equivalence. The right hand side is discrete, so we must show that $\rho^{-1}(f)$ is contractible for every morphism $f \in \Hom_{\mc{C}_{\le j}}(c_1, c_2)$. 
	
	Theorem~\ref{neckthm} gives a homotopy equivalence $\Hom_{(\mc{C}_{\le j}^\circ)_{\infty}}(c_1, c_2) \simeq \ms{Nec}_{\mc{C}_{\le j}^\circ}(c_1, c_2)$. Under this homotopy equivalence, $\rho^{-1}(f)$ corresponds to the full subcategory 
	\[
		\nec^f_{\mc{C}_{\le j}^\circ}(c_1, c_2) \subset \ms{Nec}_{\mc{C}_{\le j}^\circ}(c_1, c_2)
	\]
	characterized by the requirement that the $\mc{C}_{\le j}$-composition of the arrows in the necklace equals $f$. According to~\ref{neckreform}, $\nec^f_{\mc{C}_{\le j}^\circ}(c_1, c_2)$ can be reformulated as a category of triples $(\Delta^n, \phi, P)$ such that the $\mc{C}_{\le j}$-composition of the maps $\phi : \spine(\Delta^n) \to \mc{C}_{\le j}^\circ$ equals $f$. 
	
	\subsubsection{Forbidden sequences} 
	Suppose we are given a pair $(\Delta^n, \phi)$ as above. 
	\begin{itemize}
		\item An interval $[a, b] \subseteq [n]$ is \emph{forbidden} if $r(\phi([a, b])) < j$ and $[a, b]$ is maximal with respect to this property. 
		\item \emph{Anti-forbidden} intervals are defined similarly, using the property $r(\phi([a, b])) = j$. 
	\end{itemize}
	The condition for a triple $(\Delta^n, \phi, P)$ to be an object of $\ms{Nec}_{\mc{C}_{\le j}^\circ}(c_1, c_2)$ can be reformulated more concretely as follows: 
	\begin{itemize}
		\item $P$ intersects every forbidden interval. 
	\end{itemize}
	
	\subsubsection{Cases} 
	To show that $\ms{Nec}^f_{\mc{C}_{\le j}^\circ}(c_1, c_2)$ is contractible, we split into cases: 
	\begin{enumerate}[label=(\roman*)]
		\item $r(c_1) = r(c_2) = j$ and $f$ is basic. 
		\item $r(c_1) = r(c_2) = j$ and $f$ is not basic. 
		\item $r(c_1) < j$ and $r(c_2) = j$. 
		\item $r(c_1) = j$ and $r(c_2) < j$. 
		\item $r(c_1) < j$ and $r(c_2) < j$. 
	\end{enumerate}
	
	Assume that (i) holds. Since $f$ is basic, every sequence $\phi : \spine(\Delta^n) \to \mc{C}_{\le j}^\circ$, whose $\mc{C}_{\le j}$-composition equals $f$, must lie in $\mc{C}_j$. This implies that 
	\[
		\ms{Nec}^f_{\mc{C}_{\le j}^\circ}(c_1, c_2) = \ms{Nec}^f_{\mc{C}_j}(c_1, c_2). 
	\]
	Since $\mc{C}_j$ is a category, the right hand side is contractible, as desired. 
	
	Next, assume that (ii) holds. 
	
	\subsubsection{Removal of anti-forbidden intervals} 
	Since $f$ is not basic, every sequence $\phi : \spine(\Delta^n) \to \mc{C}_{\le j}^\circ$, whose $\mc{C}_{\le j}$-composition equals $f$, must not be contained in $\mc{C}_j$. Thus, there is at least one forbidden interval, and there are at least two anti-forbidden intervals. 
	
	Define the full subcategory 
	\[
		\ms{Nec}^{f, \ms{del}}_{\mc{C}_{\le j}^\circ}(c_1, c_2) \subset \ms{Nec}^f_{\mc{C}_{\le j}^\circ}(c_1, c_2)
	\]
	to consist of triples $(\Delta^n, \phi, P)$ such that $P$ is disjoint from the anti-forbidden intervals. This embedding is a right adjoint, and the left adjoint sends $(\Delta^n, \phi, P)$ to $(\Delta^n, \phi, P')$ where $P'$ is obtained by subtracting all anti-forbidden interval from $P$. (This does not affect the requirement that $P$ intersects every forbidden interval, because the forbidden and anti-forbidden intervals are disjoint.) Thus, the embedding is a homotopy equivalence. 
	
	Next, define the full subcategory 
	\[
		\ms{Nec}^{f, \ms{min}}_{\mc{C}_{\le j}^\circ}(c_1, c_2) \subset \ms{Nec}^{f, \ms{del}}_{\mc{C}_{\le j}^\circ}(c_1, c_2)
	\]
	to consist of triples $(\Delta^n, \phi, P)$ such that $r(\phi([1, n-1])) < j$. This implies that $[1, n-1]$ is the only forbidden interval, and the two anti-forbidden intervals are $\{0\}, \{n\}$. This embedding is a left adjoint, and the right adjoint sends $(\Delta^n, \phi, P)$ to $(\Delta^I, \phi', P)$, defined as follows: 
	\begin{itemize}
		\item $I \subset [n]$ is obtained by subtracting all anti-forbidden intervals. 
		
		Thus, the old subset $P \subseteq [n]$ lies in $I$, and $\Delta^I$ is a sub-simplex of $\Delta^n$. 
		\item For any consecutive elements $i, j \in I$, the map $\phi' : \spine(\Delta^I) \to \mc{C}_{\le j}^\circ$ sends the arrow $\Delta^{\{i, i'\}}$ to the composite of the arrows $\spine(\Delta^{[i, i']}) \xra{\phi} \mc{C}_{\le j}^\circ$. 
	\end{itemize}
	Thus, the embedding is a homotopy equivalence.

	\subsubsection{Removal of \texorpdfstring{$P$}{P}}
	Define a category $\ms{Nop}^{f, \ms{min}}$ by removing $P$ from the definition of $\ms{Nec}^{f, \ms{min}}_{\mc{C}_{\le j}^\circ}(c_1, c_2)$. In more detail: 
	\begin{itemize}
		\item An object is a pair $(\Delta^n, \phi)$, such that 
		\begin{itemize}
			\item The composite of the arrows $\phi : \spine (\Delta^n) \to \mc{C}_{\le j}^\circ$ equals $f$. 
			\item $r(\phi([1, n-1])) < j$. 
		\end{itemize}
		\item A morphism $(\Delta^{n_1}, \phi_1) \to (\Delta^{n_2}, \phi_2)$ is a map $F : \Delta^{n_1} \to \Delta^{n_2}$ such that, for every $i \in [1, n_1]$, the composite of the arrows $\spine(\Delta^{[F(i-1), F(i)]}) \xra{\phi_2} S$ equals  $\Delta^{\{i-1, i\}} \xra{\phi_1} S$.
	\end{itemize}
	\begin{lem*}
		The forgetful functor $\ms{Nec}^{f, \ms{min}}_{\mc{C}_{\le j}^\circ}(c_1, c_2) \xra{\ms{oblv}} \ms{Nop}^{f, \ms{min}}$ is a homotopy equivalence. 
	\end{lem*}
	\begin{proof}
		Observe that $\ms{oblv}$  is cocartesian, where the cocartesian arrow associated to $(\Delta^{n_1}, \phi_1, P_1)$ and $(\Delta^{n_1}, \phi_1) \xra{F} (\Delta^{n_2}, \phi_2)$ is 
		\[
			(\Delta^{n_1}, \phi_1, P_1) \xra{F} (\Delta^{n_2}, \phi_2, F(P_1)). 
		\]
		Thus, it suffices to show that each fiber $\ms{oblv}^{-1}(\Delta^n, \phi)$ is contractible. The fiber has an initial object, given by the choice $P := [1, n-1]$. 
	\end{proof}

	In fact, $\ms{Nop}^{f, \ms{min}}$ is equivalent to the category of simplices $\mb{\Delta}_{/\ms{Fact}_f}$, as follows: a simplex $\Delta^n \to \ms{Fact}_f$ is a commutative diagram 
	\begin{cd}[row sep = 4mm] 
		& d_0 \ar[d] \ar[rdd] \\
		& d_1 \ar[d] \ar[rd] \\
		c_1 \ar[ruu] \ar[ru] \ar[r] \ar[rd] \ar[rdd] & d_2 \ar[r] \ar[d] & c_2 \\
		& \vdots \ar[ru] \ar[d] \\
		& d_n \ar[ruu]
	\end{cd}
	and this corresponds to $(\Delta^{n+2}, \phi) \in \ms{Nop}^{f, \ms{min}}$, where $\phi$ is given by 
	\[
	c_1 \to d_0 \to d_1 \to d_2 \to \cdots \to d_n \to c_2.
	\]
	According to~\cite[Prop.\ 7.3.15]{cis} and \cite[Prop.\ 7.1.10]{cis}, $\mb{\Delta}_{/\ms{Fact}_f}$ is homotopy equivalent to $\ms{Fact}_f$, which is contractible by hypothesis. Thus, $\ms{Nec}^f_{\mc{C}_{\le j}^\circ}(c_1, c_2)$ is contractible, so case (ii) is done. 
	
	\subsubsection{Proof for cases (iii), (iv), (v)} 
	As before, we show that $\ms{Nec}^f_{\mc{C}_{\le j}^\circ}(c_1, c_2)$ is homotopy equivalent to a category $\ms{Nop}^{f, \ms{min}}$ whose objects are pairs $(\Delta^n, \phi)$ where $r(\phi([1, n-1])) < j$. This time, $\ms{Nop}^{f, \ms{min}}$ is contractible because it has an initial object, given by $(\Delta^1, f)$. This concludes the proof of Theorem~\ref{recognize}.

	\mysubsec{\texorpdfstring{$\infty$}{oo}-bistratified colimits} \label{extra-bis} 
	As in the special case of Reedy categories, an essential feature of $\infty$-bistratified categories is that colimits indexed by them can be computed inductively. We will formulate a special case of this principle which is convenient for applications. 
	
	Work in the setting of~\ref{ss-bis}, and fix $j \in J$. Suppose we are given a functor $F : \mc{C}_{\le j} \to \mc{E}$ to an $\infty$-category which admits colimits. 
	
	\subsubsection{Relative functor}  \label{bis-relative} 
	Define the functor $\bar{F} : \mc{C}_{\le j} \to \mc{E}$ via the following pushout square: 
	\begin{cd}
		\LKE_{\mc{C}_{<j} \hra \mc{C}_{\le j}} \Res_{\mc{C}_{<j} \hra \mc{C}_{\le j}} F \ar[r] \ar[d] & F \ar[d] \\
		\Res_{\mc{C}_{\le j} \to \pt} {\displaystyle\colim_{\mc{C}_{<j}}} F \ar[r, "\varphi"] & \bar{F}
	\end{cd} 
	The bottom-left vertex is the constant functor with value $\colim_{\mc{C}_{<j}} F$. The value of this square at $c \in  \mc{C}_{\le j}$ is 
	\begin{cd}
		\displaystyle\colim_{\la \mc{C}_{<j} \xra{\mc{C}} c \ra} F \circ \ms{tail} \ar[r] \ar[d] & F(c) \ar[d] \\
		\displaystyle\colim_{\mc{C}_{<j}} F \ar[r, "\varphi(c)"] & \bar{F}(c)
	\end{cd} 	
	Let $\tilde{\mc{E}} := \mc{E}_{(\colim_{\mc{C}_{<j}} F)/}$, and denote its initial object by $\hat 0$. We can view $(\bar{F}, \varphi)$ as a functor $\tilde{F} : \mc{C}_{\le j} \to \tilde{\mc{E}}$. Since $\mc{C}_{<j} \to \mc{C}_{\le j}$ is fully faithful, $\tilde{F}$ sends $\mc{C}_{<j}$ to $\hat 0$. 
	
	\begin{lem*}
		There is an $\tilde{\mc{E}}$-isomorphism $\colim_{\mc{C}_{\le j}} \tilde{F} \simeq \Big(\colim_{\mc{C}_{<j}} F \to \colim_{\mc{C}_{\le j}} F\Big)$. 
	\end{lem*}
	\begin{proof}
		The colimit of the previous square is a pushout square in $\mc{E}$
		\begin{cd}
			\displaystyle\colim_{\mc{C}_{<j}} F \ar[r] \ar[d] & \displaystyle\colim_{\mc{C}_{\le j}} F \ar[d] \\
			{|\mc{C}_{\le j}|} \otimes \displaystyle\colim_{\mc{C}_{<j}} F \ar[r] & \displaystyle\colim_{\mc{C}_{\le j}} \bar{F}
		\end{cd}
		where $|\mc{C}_{\le j}| \in \infty\text{-}\ms{Grpd}$ is obtained from $\mc{C}_{\le j}$ by inverting all morphisms. If we write 
		\[
			\colim_{\mc{C}_{\le j}} \tilde{F} = \Big( \colim_{\mc{C}_{<j}} F \xra{i} y \Big) \in \tilde{\mc{E}}, 
		\]
		then the computation of colimits in $\tilde{\mc{E}}$ implies that the following is a pushout square in $\mc{E}$ 
		\begin{cd}
			{|\mc{C}_{\le j}|} \otimes \displaystyle\colim_{\mc{C}_{<j}} F \ar[r] \ar[d] & \displaystyle\colim_{\mc{C}_{\le j}} \bar{F} \ar[d] \\
			\colim_{\mc{C}_{<j}} F \ar[r, "i"] & y
		\end{cd}
		Stacking these two squares gives the desired result. 
	\end{proof}
	
	Thus, computing $\colim_{\mc{C}_{\le j}} F$ in terms of $\colim_{\mc{C}_{<j}} F$ is equivalent to computing $\colim_{\mc{C}_{\le j}} \tilde{F}$. 
	
	\subsubsection{Downward subcategory} \label{bis-downward} 
	Let $\mc{C}_{\le j}^{\downarrow} \subset \mc{C}_{\le j}$ be the subcategory consisting of all objects and all arrows except for the non-basic level-$j$ arrows and all arrows from $\mc{C}_{<j}$ to $\mc{C}_j$. In other words, every arrow in $\mc{C}_{\le j}^{\downarrow}$ satisfies exactly one of the following: it lies in $\mc{C}_{<j}$, it lies in $\mc{C}_j$, or it goes from $\mc{C}_j$ to $\mc{C}_{<j}$. 
	
	\begin{lem*}
		We have $\colim_{\mc{C}_{\le j}} \tilde{F} \simeq \colim_{\mc{C}_{\le j}^{\downarrow}} \tilde{F}$. 
	\end{lem*}
	\begin{proof}
		The embedding $\mc{C}_{\le j}^{\downarrow} \hra \mc{C}_{\le j}$ factors as 
		\begin{cd}
			\mc{C}_{\le j}^{\downarrow} \ar[r, "\iota_1"] & \mc{C}_{\le j}^\circ \ar[r, "\iota_2"] & \mc{C}_{\le j} 
		\end{cd}
		This induces the following two maps, which we claim are isomorphisms. 
		\begin{cd}
			\displaystyle\colim_{\mc{C}_{\le j}^{\downarrow}} \tilde{F} \ar[r, "\jmath_1"] & \displaystyle\colim_{\mc{C}_{\le j}^\circ} \tilde{F} \ar[r, "\jmath_2"] & \displaystyle\colim_{\mc{C}_{\le j}} \tilde{F}
		\end{cd}
		The map $\jmath_2$ is an isomorphism because $\iota_2$ is a categorical equivalence. The map $\jmath_1$ is an isomorphism because every simplex $\sigma : \Delta^m \to \mc{C}_{\le j}^\circ$ which does not factor through $\iota_1$ satisfies $\tilde{F}(\sigma(0)) = \hat0$, so this simplex does not affect the colimit. 
	\end{proof}
	
	\subsubsectiona 
	The advantage in working with $\colim_{\mc{C}_{\le j}^{\downarrow}}$ is that it admits a map to $\Delta^1$, where an object maps to 0 if and only if its rank equals $j$. The next result pertains to this situation. 
	
	\begin{prop*} \label{constprop} 
		Let $F: \mc{B} \to \mc{D}$ be a functor between categories, where $\mc{D}$ admits colimits. Suppose we are given a map $\mc{B} \to \Delta^1$, whose fibers are denoted $\mc{B}_0$ and $\mc{B}_1$. Then there is a pushout square in $\mc{D}$ 
		\begin{cd}
			\displaystyle\colim_{\la \mc{B}_0 \xra{\mc{B}} \mc{B}_1 \ra} F \circ \ms{tail} \ar[r] \ar[d] & \displaystyle\colim_{\mc{B}_1} F \ar[d] \\
			\displaystyle\colim_{\mc{B}_0} F \ar[r] & \displaystyle\colim_{\mc{B}} F
		\end{cd}
	\end{prop*}
	\begin{proof}
		Consider the opposite twisted arrow category $\ms{Tw}^{\ms{op}}(\mc{B})$, in which a morphism $f \to g$ corresponds to a commutative diagram in $\mc{B}$ 
		\begin{cd}
			\bullet \ar[r] \ar[d, "f"] & \bullet \ar[d, "g"] \\
			\bullet \ar[r, leftarrow] & \bullet 
		\end{cd}
		The functor $\ms{tail} : \ms{Tw}^{\ms{op}}(\mc{B}) \to \mc{B}$ is a cartesian fibration with contractible fibers, so it is final. Thus, we have 
		\[
		\colim_{\ms{Tw}^{\ms{op}}(\mc{B})} F \circ \ms{tail} \simeq \colim_{\mc{B}} F. 
		\]
		
		Define the functor $\pi : \ms{Tw}^{\ms{op}}(\mc{B}) \to \Lambda^2_0 := (1 \leftarrow 0 \rightarrow 2)$ as follows:
		\[
		\pi(b \xra{f} b') = \begin{cases}
			0 & \text{if $b \in \mc{B}_0$ and $b' \in \mc{B}_1$} \\
			1 & \text{if $b \in \mc{B}_1$ and $b' \in \mc{B}_1$} \\
			2 & \text{if $b \in \mc{B}_0$ and $b' \in \mc{B}_0$}	
		\end{cases}
		\]
		If we define $G := \LKE_\pi (F \circ \ms{tail})$, then $\colim_{\Lambda^2_0} G \simeq \colim_{\ms{Tw}^{\ms{op}}(\mc{B})} F \circ \ms{tail}$, so we have a pushout square in $\mc{D}$ 
		\begin{cd}
			G(0) \ar[r] \ar[d] & G(1) \ar[d] \\
			G(2) \ar[r] & \displaystyle\colim_{\mc{B}} F
		\end{cd}
		It remains to compute the values of $G$, which are given as follows: 
		\e{
			G(0) &\simeq \colim_{\pi^{-1}(0)} F \circ \ms{tail} \\
			G(1) &\simeq \colim_{\pi^{-1}(\{0, 1\})} F \circ \ms{tail} \\
			G(2) &\simeq \colim_{\pi^{-1}(\{0, 2\})} F \circ \ms{tail}. 
		} 
		
		Since the cartesian fibrations $\pi^{-1}(0) \xra{\ms{tail}} \mc{B}_0$ and $\la \mc{B}_0 \xra{\mc{B}} \mc{B}_1\ra \xra{\ms{tail}} \mc{B}_0$ correspond to isomorphic functors $\mc{B}_0 \to \infty\text{-}\ms{Grpd}$, the colimit for $G(0)$ rewrites as follows: 
		\[
			G(0) \simeq \colim_{\la \mc{B}_0 \xra{\mc{B}} \mc{B}_1\ra} F \circ \ms{tail}. 
		\]
		
		We claim that the embedding $\pi^{-1}(1) \hra \pi^{-1}(\{0, 1\})$ is final. For this, we fix $(b_0 \xra{f} b_1) \in \pi^{-1}(\{0, 1\})$ and show that $\la (b_0 \xra{f} b_1) \xra{\pi^{-1}(\{0, 1\})} \pi^{-1}(1) \ra$ is contractible. An object of this undercategory is a diagram 
		\begin{cd}
			b_0 \ar[r] \ar[d, "f"] & b'_0 \ar[d, "g"] \\
			b_1 \ar[r, leftarrow] & b'_1
		\end{cd}
		Consider the full subcategory defined by the requirement that $b'_1 = b_1$. Its embedding is a left adjoint. It also has a terminal object ($b'_0 = b'_1 = b_1$), so the contractibility follows.  
		
		The claim implies $G(1) \simeq \colim_{\pi^{-1}(1)} F \circ \ms{tail}$. Since $\pi^{-1}(1) = \ms{Tw}^{\ms{op}}(\mc{B}_1)$, the previous reasoning implies that this is equivalent to $\colim_{\mc{B}_1} F$. Similarly, one can show that $G(2) \simeq \colim_{\mc{B}_2} F$, as desired. 
	\end{proof}
	
	\begin{rmk*}
		The proposition also works for $\infty$-categories, but we will not need this. 
	\end{rmk*}
	
	\subsubsection{Two important cases}  \label{cases} 
	We will only apply the above results in these two cases: 
	\begin{enumerate}[label=(\roman*)]
		\item $\tilde{F}$ sends $\mc{C}_j$-maps to isomorphisms. Also, $\la \mc{C}_j \xra{\mc{C}} \mc{C}_{<j} \ra$ is empty, and 1$\mc{C}_j$ is contractible. 
		\item $\tilde{F}$ sends $\mc{C}_j$-maps to isomorphisms. The categories $\mc{C}^\downarrow_{\le j}$ and $\mc{C}_{<j}$ are contractible. 
	\end{enumerate}
	In case (i), the second hypothesis implies $\colim_{\mc{C}^\downarrow_{\le j}} \tilde{F} \simeq \colim_{\mc{C}_j} \tilde{F}$. The first hypothesis implies that, for any $c \in \mc{C}_j$, we have $\tilde{F}(c) \simeq \colim_{\mc{C}_j} \tilde{F}$, see~\cite[Cor.\ 4.4.4.10]{htt}. Now~\ref{bis-relative} gives the following pushout square in~$\mc{E}$. 
	\begin{cd}
		\displaystyle\colim_{\la \mc{C}_{<j} \xra{\mc{C}} c \ra} F \circ \ms{tail} \ar[r] \ar[d] & F(c) \ar[d] \\
		\displaystyle\colim_{\mc{C}_{< j}} F \ar[r] & \displaystyle\colim_{\mc{C}_{\le j}} F
	\end{cd}
	
	In case (ii), the second hypothesis implies that the map $\ms{tail} : \la \mc{C}_j \xra{\mc{C}} \mc{C}_{<j} \ra \to \mc{C}_j$ is a homotopy equivalence. (To see this, apply Proposition~\ref{constprop} to the constant functor with value $\pt \in \ms{Spaces}$.) Then the first hypothesis implies that $\colim_{\mc{C}^\downarrow_{\le j}} \tilde{F} \simeq \colim_{\mc{C}_{<j}} \tilde{F}$ via a second application of Proposition~\ref{constprop}. Since this is $\hat0$, we conclude that 
	\[
		\colim_{\mc{C}_{<j}} F \simeq \colim_{\mc{C}_{\le j}} F. 
	\]

	\mysubsec{Rigidly bistratified categories} \label{rigidcat} 
	
	\subsubsectiona  \label{rigidcat-gen} 
	Here is a more obvious generalization of the notion of a Reedy category. It turns out to be a special case of~\ref{ss-bis}. 
	
	A \emph{generating datum} consists of the following: 
	\begin{itemize}
		\item A ground set $X$. 
		\item Two categories $\Domc$ and $\Embc$ with $\ms{Obj}(\Domc) = \ms{Obj}(\Embc) = X$. 
		
		Arrows in $\Embc$ are drawn as $\rightarrowtail$, and arrows in $\Domc$ are drawn as $\twoheadrightarrow$. 
		\item \emph{(Exchange map.)} For each diagram on the left, a choice of diagram on the right: 
		\[
		\begin{tikzcd}
		x_1 \ar[r, rightarrowtail] & x_2 \ar[d, twoheadrightarrow] \\
		& x_3
		\end{tikzcd} \quad \rightsquigarrow \quad 
		\begin{tikzcd}
		x_1 \ar[d, twoheadrightarrow] \\
		x_4 \ar[r, rightarrowtail] & x_3
		\end{tikzcd}
		\]
	\end{itemize}
	We impose the following axioms: 
	\begin{itemize}
		\item If the following two solid diagrams are equal, then the outer diagrams are equal: 
		\[
		\begin{tikzcd}
		& x_2 \ar[rd, rightarrowtail] \ar[d, dashed, twoheadrightarrow] \\
		x_1 \ar[ru, rightarrowtail] \ar[d, dashed, twoheadrightarrow] & \bullet  \ar[rd, dashed, rightarrowtail] & x_3 \ar[d, twoheadrightarrow] \\
		\bullet \ar[ru, dashed, rightarrowtail] \ar[rr, dashed, rightarrowtail] & & x_4
		\end{tikzcd} 
		\qquad \qquad
		\begin{tikzcd}
		& x_2 \ar[rd, rightarrowtail] \\
		x_1 \ar[ru, rightarrowtail] \ar[rr, dashed, rightarrowtail] \ar[d, dashed, twoheadrightarrow] && x_3 \ar[d, twoheadrightarrow] \\
		\bullet \ar[rr, dashed, rightarrowtail] & & x_4
		\end{tikzcd} 
		\]
		The squares come from the exchange map, and the triangles come from composition in $\Embc$. Also, we impose the analogous axiom for one $\rightarrowtail$ and two $\twoheadrightarrow$'s. 
		\item The exchange map interacts with identity morphisms in the obvious way. 
	\end{itemize}
	We generate a category $\Ic$ as follows: the object set is $X$, the generating arrows are given by $\Domc$ and $\Embc$, and relations between arrows are given by the exchange map. Thanks to the axioms, this category can be concretely described as follows: a morphism is a diagram $x_1 \sra \bullet \rightarrowtail x_2$, and composition of morphisms is given by 
	\begin{cd}
		x_1 \ar[d, twoheadrightarrow] \\
		\bullet \ar[r, rightarrowtail] \ar[d, dashed, twoheadrightarrow] & x_2 \ar[d, twoheadrightarrow] \\
		\bullet \ar[r, dashed, rightarrowtail] & \bullet \ar[r, rightarrowtail] & x_3
	\end{cd}
	By construction, every $\ms{Rig}^{\ms{c}}$-morphism factors uniquely and functorially as a $\ms{Dom}^{\ms{c}}$-morphism followed by an $\ms{Emb}^{\ms{c}}$-morphism. 
	
	\begin{lem*}
		Suppose we are given a map $r : X \to J$ to an ordered set, such that 
		\begin{itemize}
			\item $r$ is strictly increasing along non-identity morphisms in $\Embc$. 
			\item $r$ is weakly decreasing along morphisms in $\Domc$. 
		\end{itemize} 
		Then $\ms{Rig}^{\ms{c}}$ is $\infty$-bistratified with rank function $r$. 
	\end{lem*}
	\begin{proof}
		Let $f := x_1 \sra x \rightarrowtail x_2$ be a level $\ms{Rig}^{\ms{c}}$-morphism, with $r(x_1) = r(x_2) = j$. Then $r(x) \le j$. We claim that $f$ is basic if and only if $r(x) = j$. Clearly, $r(x) < j$ implies that the morphism is not basic. Conversely, suppose that $f$ is not basic, so that it factors as 
		\[
			x_1 \sra y_1 \rightarrowtail y_2 \sra y_3 \rightarrowtail x_2
		\]
		where $r(y_2) < j$. This inequality implies $r(y_1), r(y_3) < j$. Using the exchange map, we can write this as 
		\[
			x_1 \sra y_1 \sra y'_2 \rightarrowtail y_3 \rightarrowtail x_2
		\]
		where $r(y'_2) < j$. By uniqueness, the factorization $x_1 \sra y'_2 \rightarrowtail x_2$ must equal $x_1 \sra x \rightarrow x_2$, so $r(x) = r(y'_2) < j$, as desired. 
		
		Note that $r(x) = j$ and $r(x_2) = j$ implies that $x \rightarrowtail x_2$ is an identity. Thus, the set of basic level morphisms equals the set of level morphisms in $\ms{Dom}^{\ms{c}}$. Since this is a subcategory of $\ms{Rig}^{\ms{c}}$, we have proven~\ref{cond-b1}. 
		
		Next, we prove~\ref{cond-b2}. Assume that $f$ is not basic. We want to show that $\ms{Fact}_f$ is contractible. An object in $\ms{Fact}_f$ may be denoted as $x_1 \to y_2 \to x_2$ (see above), where the first morphism factors as $x_1 \sra y_1 \rightarrowtail y_2$, and the second morphism factors as $y_2 \sra y_3 \rightarrowtail x_2$. Consider the full subcategory 
		\[
			\ms{Fact}_f^{\ms{post}} \subset \ms{Fact}_f
		\]
		which is defined by the requirement that $y_2 \to x_2$ lies in $\ms{Emb}^{\ms{c}}$. The embedding is a right adjoint, and the left adjoint sends the aforementioned factorization to $x_1 \to y_3 \rightarrowtail x_2$. Furthermore, $\ms{Fact}_f^{\ms{post}}$ has an initial object given by the canonical factorization $x_1 \sra x \rightarrowtail x_2$. This proves the desired contractibility. 
	\end{proof}

	\begin{ex*}
		The category of finite type schemes over a field $k$ can be presented in this way. The subcategory $\ms{Emb}^{\ms{c}}$ consists of closed embeddings, and $\ms{Dom}^{\ms{c}}$ consists of maps such that the scheme-theoretic image equals the target. The latter are \emph{dominant} maps of schemes, which motivates our notation. 
	\end{ex*}

	\subsubsectiona  \label{rigidcat-deform} 
	Suppose we are given a generating datum for which $\ms{Emb}^{\ms{c}}$ is right-cancellable, i.e.\ $f_1 \circ f = f_2 \circ f$ implies $f_1 = f_2$. In this case, we will define what it means to `deform' the generating datum. 
	
	\begin{rmk*}
		For the simplex category, $\ms{Emb}^{\ms{c}}$ consists of strictly increasing maps between finite nonempty ordered sets, so it does not satisfy right-cancellability. For schemes, $\ms{Emb}^{\ms{c}}$ consists of closed embeddings, and it also does not satisfy right-cancellability. 
	\end{rmk*}
	
	A \emph{deformation} of the generating datum consists of a subcategory $\Embd \subset \Embc$ and a \emph{turning map} which sends arrows in $\ms{Dom}^{\ms{c}}$ to arrows in $\Embd$ as shown: 
	\[
	\begin{tikzcd}
	x_1 \ar[d, twoheadrightarrow] \\
	x_2
	\end{tikzcd} \quad \rightsquigarrow \quad 
	\begin{tikzcd}
	\phantom{a} \\
	x_2 \ar[r, hookrightarrow, "\mr{turn}"] & x_3
	\end{tikzcd}
	\]
	Arrows in $\Embd$ will be drawn as $\hookrightarrow$. We impose the following axioms: 
	\begin{itemize}
		\item Suppose we are given a solid diagram as below: 
		\begin{cd}
			x_1 \ar[d, twoheadrightarrow] \\
			x_2 \ar[d, dashed, twoheadrightarrow] \ar[r, dashed, hookrightarrow, "\mr{turn}"] & x_3 \ar[d, twoheadrightarrow] \\
			\bullet \ar[r, dashed, rightarrowtail] & x_4 \ar[r, dashed, hookrightarrow, "\mr{turn}"] & x_5
		\end{cd}
		Fill in the square using the exchange map. Then we require that the turn of the vertical composite equals the horizontal composite. 
		\item Suppose we are given a solid diagram as below: 
		\begin{cd}
			& x_2 \ar[d, twoheadrightarrow] \\
			x_1 \ar[ru, hookrightarrow] \ar[d, dashed, twoheadrightarrow] & x_3 \ar[r, hookrightarrow, dashed, "\mr{turn}"] & \bullet \\
			\bullet \ar[ru, dashed, rightarrowtail, start anchor = {[yshift = 1ex]}] \ar[r, dashed, hookrightarrow, swap, "\mr{turn}"] & \bullet \ar[ru, dotted, hookrightarrow, swap, "\exists"]
		\end{cd}
		Fill in the square using the exchange map, and create the horizontal $\hookrightarrow$'s using the turning map. We require that there exists a dotted arrow in $\Embd$ making the lower square commute. (It is unique by right-cancellability.) 
	\end{itemize} 
	
	The axioms imply that we can define a new generating datum on the same ground set $X$, with new categories $\Embd$ and $\Domd$ given as follows: 
	\begin{itemize}
		\item $\Domd$ is the subcategory of $\ms{Rig}^{\ms{c}}$ consisting of the maps 
		\begin{cd}
			x_1 \ar[d, twoheadrightarrow] \\
			\bullet \ar[r, dashed, hookrightarrow, "\mr{turn}"] & x_2
		\end{cd}
		The first axiom says that $\Domd$ is closed under composition. 
		\item The exchange map is defined using the second axiom. 
	\end{itemize}
	We leave it as an exercise to show that this new generating datum satisfies~\ref{rigidcat-gen}. 
	
	The new generating datum gives a new category $\ms{Rig}^{\ms{d}}$. We think of it as a deformation of $\ms{Rig}^{\ms{c}}$ because the $\ms{Dom}^{\ms{d}}$-maps are `rotated' versions of the $\ms{Dom}^{\ms{c}}$-maps. (The `rotation' is specified by the turning map.) In fact, $\ms{Rig}^{\ms{d}}$ identifies with the subcategory of $\ms{Rig}^{\ms{c}}$ consisting of morphisms $x_1 \sra x \rightarrowtail x_2$ such that $x \rightarrowtail x_2$ factors as $x \overset{\mr{turn}}{\hra} \bullet \hra x_2$, where the `turn' arrow comes from $x_1 \sra x$. (The factorization is unique by right-cancellability.) 
	
	A function $r : X \to J$ to an ordered set is called \emph{good} if the following hold: 
	\begin{itemize}
		\item $r$ is strictly increasing along non-identity morphisms in $\Embd$. 
		\item $r$ is weakly decreasing along morphisms in $\Domd$. 
		\item For each $\ms{Dom}^{\ms{d}}$-map $x_1 \sra x \overset{\mr{turn}}{\hra} x_2$, we have $r(x_1) = r(x)$ if and only if $r(x_1) = r(x_2)$.  
	\end{itemize}
	The first two bullets imply that $\ms{Rig}^{\ms{d}}$ is $\infty$-bistratified with rank function $r$, by Lemma~\ref{rigidcat-gen}. Furthermore, the third bullet implies (via the proof of that lemma) that 
	\[
		(\text{basic level morphisms of } \ms{Rig}^{\ms{d}}) = (\text{level morphisms of } \ms{Dom}^{\ms{d}}) = (\text{level morphisms of } \ms{Dom}^{\ms{c}}). 
	\]

	\subsubsection{Relative functor}  \label{rigid-rel} 
	Fix a good rank function $r : X \to J$ as above. Suppose we are given a functor $F : \ms{Rig}^{\ms{d}} \to \mc{E}$ to an $\infty$-category which admits colimits. As in~\ref{extra-bis}, we fix $j \in J$ and consider the relative functor $\tilde{F} : \ms{Rig}^{\ms{d}, \downarrow}_{\le j} \to \tilde{\mc{E}}$. Note that $\ms{Dom}^{\ms{d}}_{\le j} \subset \ms{Rig}^{\ms{d}, \downarrow}_{\le j}$. 
	
	\begin{lem*}
		We have $\colim_{\mc{C}^\downarrow_{\le j}} \tilde{F} \simeq \colim_{\ms{Dom}^{\ms{d}}_{\le j}} \tilde{F}$. 
	\end{lem*}
	\begin{proof}
		Applying Proposition~\ref{constprop} to $\ms{Rig}^{\ms{d}, \downarrow}_{\le j}$ and $\ms{Dom}^{\ms{d}}_{\le j}$ gives pushout squares in $\tilde{\mc{E}}$ as follows: 
		\[	
		\begin{tikzcd}
			\displaystyle\colim_{\la \ms{Rig}^{\ms{d}}_j \xra{\ms{Rig}^{\ms{d}}} \ms{Rig}^{\ms{d}}_{<j} \ra} \tilde{F} \circ \ms{tail} \ar[r] \ar[d] & \displaystyle\colim_{\ms{Rig}^{\ms{d}}_{<j}} \tilde{F} \ar[d] \\
			\displaystyle\colim_{\ms{Rig}^{\ms{d}}_{j}} \tilde{F} \ar[r] & \displaystyle\colim_{\ms{Rig}^{\ms{d}}_{\le j}} \tilde{F}
		\end{tikzcd}
		\qquad  
		\begin{tikzcd}
		\displaystyle\colim_{\la \ms{Dom}^{\ms{d}}_j \xra{\ms{Dom}^{\ms{d}}} \ms{Dom}^{\ms{d}}_{<j} \ra} \tilde{F} \circ \ms{tail} \ar[r] \ar[d] & \displaystyle\colim_{\ms{Dom}^{\ms{d}}_{<j}} \tilde{F} \ar[d] \\
		\displaystyle\colim_{\ms{Dom}^{\ms{d}}_{j}} \tilde{F} \ar[r] & \displaystyle\colim_{\ms{Dom}^{\ms{d}}_{\le j}} \tilde{F}
		\end{tikzcd}
		\]
		The upper-right colimits are $\hat0$ because $\tilde{F}$ sends these categories to $\hat0$. The lower-left colimits coincide because $\ms{Rig}^{\ms{d}}_{j} = \ms{Dom}^{\ms{d}}_{j}$. Thus, it suffices to show that the upper-left colimits coincide. 
		
		The full embedding $\la\Domd_j \xra{\Domd} \Domd_{<j}\ra \hra \la\Id_j \xra{\Id} \Id_{<j}\ra$ is a left adjoint, with right adjoint given by the functorial factorization of each $\ms{Rig}^{\ms{d}}$-map into a $\ms{Dom}^{\ms{d}}$-map followed by an $\ms{Emb}^{\ms{d}}$-map. The functor $\tilde{F} \circ \ms{tail}$ sends the counit of this adjunction to a natural isomorphism, because the aforementioned factorization does not change the tail of the arrow. Now Lemma~\ref{adjoint-colim} (below) implies that the upper-left colimits coincide. 
	\end{proof}
	
	\begin{lem} \label{adjoint-colim} 
		Suppose we are given an adjunction $L : \mc{C} \rightleftarrows \mc{D} : R$ and a functor $F: \mc{D} \to \mc{E}$ to an $\infty$-category which admits colimits. If $F$ sends the counit $LR \Rightarrow \id_{\mc{D}}$ to a natural isomorphism, then the map $\colim FL \to \colim F$ is an isomorphism. 
	\end{lem}
	\begin{proof}
		The hypothesis implies that $\LKE_L \Res_L F \Rightarrow F$ is a natural isomorphism. The result follows from left Kan extension along $\mc{D} \to \pt$. 
	\end{proof}
	
	\subsubsection{Anti-deformation of the relative functor}  \label{rigid-anti} 
	Assume that $\ms{Emb}^{\ms{c}}$ is a groupoid. Then, for every solid diagram (as below), there exists a unique choice of dashed arrows such that the square comes from the exchange map: 
	\begin{cd}
		\bullet \ar[r, rightarrowtail] \ar[d, twoheadrightarrow] & \bullet \ar[d, dashed, twoheadrightarrow] \\
		\bullet \ar[r, dashed, rightarrowtail] & \bullet
	\end{cd}
	The idea is to invert the given $\rightarrowtail$ arrow and apply the exchange map. 
	
	Next, using the functor $\tilde{F} : \ms{Dom}^{\ms{d}}_{\le j} \to \tilde{\mc{E}}$, we define an \emph{anti-deformed} functor $\tilde{F}^{\ms{c}} : \ms{Dom}^{\ms{c}}_{\le j} \to \tilde{\mc{E}}$ as follows: 
	\begin{itemize}
		\item We may assume that the restriction of $\tilde{F}$ to $\ms{Dom}^{\ms{d}}_{< j}$ is constant with value $\hat0 \in \tilde{\mc{E}}$. 
		\item For each simplex $\sigma : \Delta^n \to \ms{Dom}^{\ms{c}}_{\le j}$, the resulting simplex $\tilde{F}^{\ms{c}} \circ \sigma : \Delta^n \to \tilde{\mc{E}}$ is defined by splitting into cases: 
		\begin{itemize}
			\item Assume that there exists $i \in [n]$ such that $r(\sigma(i)) = j$. Let $i$ be the maximal such index. Define a `deformed' simplex $\sigma^{\ms{d}} : \Delta^n \to \ms{Dom}^{\ms{d}}_{\le j}$ by modifying the $[i+1, n]$-part of $\sigma$ as follows: 
			\begin{cd}
				\sigma(i) \ar[d, twoheadrightarrow] \\
				\sigma(i+1) \ar[d, twoheadrightarrow] \ar[r, dashed, hookrightarrow, "\mr{turn}"] & \sigma^{\ms{d}}(i+1) \ar[d, dashed, twoheadrightarrow] \\
				\sigma(i+2) \ar[d, twoheadrightarrow] \ar[r, dashed, rightarrowtail] & \bullet \ar[d, dashed, twoheadrightarrow] \ar[r, dashed, hookrightarrow, "\mr{turn}"] & \sigma^{\ms{d}}(i+2) \ar[d, dashed, twoheadrightarrow] \\
				\vdots & \vdots & \vdots & \ddots
			\end{cd}
			(The squares are filled in using the previous paragraph. If $i = n$, then $\sigma^{\ms{d}} = \sigma$.) Define $\tilde{F}^{\ms{c}} \circ \sigma := \tilde{F} \circ \sigma^{\ms{d}}$. 
			\item Assume that no such $i$ exists, i.e.\ $\sigma$ factors through $\ms{Dom}^{\ms{c}}_{<j}$. Define $\tilde{F}^{\ms{c}} \circ \sigma$ to be constant with value $\hat0 \in \tilde{\mc{E}}$. 
		\end{itemize}
	\end{itemize}
	It is easy to see that this defines a valid functor $\tilde{F}^{\ms{c}}$, whose restriction to $\ms{Dom}^{\ms{c}}_j = \ms{Dom}^{\ms{d}}_j$ agrees with $\tilde{F}$, and whose restriction to $\ms{Dom}^{\ms{c}}_{< j}$ is constant with value $\hat0 \in \tilde{\mc{E}}$. 
	
	\begin{lem*}
		We have $\colim_{\ms{Dom}^{\ms{d}}_{\le j}} \tilde{F} \simeq \colim_{\ms{Dom}^{\ms{c}}_{\le j}} \tilde{F}^{\ms{c}}$. 
	\end{lem*}
	\begin{proof}
		As in the previous lemma, we reduce to showing that 
		\[
			\colim_{\la \ms{Dom}^{\ms{d}}_j \xra{\ms{Dom}^{\ms{d}}} \ms{Dom}^{\ms{d}}_{<j} \ra} \tilde{F} \circ \ms{tail} \simeq \colim_{\la \ms{Dom}^{\ms{c}}_j \xra{\ms{Dom}^{\ms{c}}} \ms{Dom}^{\ms{c}}_{<j} \ra} \tilde{F}^{\ms{c}} \circ \ms{tail}.
		\]
		This follows from the equivalence $\la\Domc_j \xra{\Domc} \Domc_{<j}\ra \xra{\sim} \la\Domd_j \xra{\Domd} \Domd_{<j}\ra$, which sends a $\Domc$-map $x_1 \twoheadrightarrow x_2$ to the $\Domd$-map $x_1 \twoheadrightarrow x_2 \overset{\mr{turn}}{\hra} x_3$. 
	\end{proof}
	
	Thus, if we understand $\ms{Dom}^{\ms{c}}$, then we can compute colimits indexed by $\ms{Rig}^{\ms{d}}$.

	\section{Demazure category}  \label{dem} 
	
	\mysubsec{Main construction} 
	
	\subsubsectiona
	Fix a Coxeter system $(I, W)$, i.e.\ $W$ is a Coxeter group and $I$ indexes simple reflections. 
	\begin{itemize}
		\item Let $\ell : W \to \BZ_{\ge 0}$ be the usual length function. 
		\item Denote the (strong) Bruhat order by $\preceq$. 
		\item The usual product on $W$ is called the \emph{Coxeter} product. 
		\item There is also the \emph{Demazure} product, which makes $W$ into a monoid.  
		\item For each subset $J \subseteq I$, let $W_J \subseteq W$ be the subgroup generated by the simple reflections in $J$. If $W_J$ is finite, we say that $J$ is \emph{finite type}. 
	\end{itemize}
	
	\subsubsection{Coxeter category} We will apply~\ref{rigidcat-gen}. 
	
	Let $X$ be the set of triples $([1, n], w, J)$ where $w = (w_1, \ldots, w_n)$ is a sequence in $W$, and $J = (J_1, \ldots, J_n)$ is a sequence of finite type subsets of $I$, such that $w_i \in W_{J_i}$ for each $i$. These are called \emph{words}. 
	
	Define $\ms{Dom}^{\ms{c}}$ as follows. A $\ms{Dom}^{\ms{c}}$-morphism $([1, n], w, J) \to ([1, n'], w', J')$ is a weakly-increasing map $\psi : [1, n] \to [1, n']$ such that, for each $j \in [1, n']$, we have $w'_j = \prod_{i \in \psi^{-1}(j)} w_i$ and $J'_j \supseteq \bigcup_{i \in \psi^{-1}(j)} J_i$. If $\psi^{-1}(j)$ is empty, we declare that the empty product is $1 \in W$, and the empty union is $\emptyset \subset I$. 
	
	Define $\ms{Emb}^{\ms{c}}$ as follows. For any two words of the form $([1, n], w, J), ([1, n], w', J)$, there is exactly one $\ms{Emb}^{\ms{c}}$ morphism from the first to the second. Note that $\ms{Emb}^{\ms{c}}$ is a groupoid. 
	
	To define the exchange map, take any solid diagram 
	\begin{cd}
		([1,n^0],w^0,J^0) \ar[r, rightarrowtail] \ar[d, dashed, twoheadrightarrow, "{\psi : [1,n^0] \to [1,n^2]}"] & ([1,n^0],w^1,J^0) \ar[d, twoheadrightarrow, "{\psi : [1,n^0] \to [1,n^2]}"] \\
		([1,n^2],w^3,J^2) \ar[r, dashed, rightarrowtail] & ([1,n^2],w^2,J^2)
	\end{cd}
	and produce the dashed diagram by setting $w^3_j := \prod_{i \in \psi^{-1}(j)} w^0_i$ for every $j \in [1, n^2]$. The axiom with two $\rightarrowtail$ arrows is obvious, and the axiom with two $\twoheadrightarrow$ arrows expresses the associativity of multiplication in $W$. 
	
	Then $\ms{Rig}^{\ms{c}}$ is called the \emph{Coxeter category}, because it encodes the Coxeter product on $W$. 
	
	\subsubsection{Rank function} 
	Define the subset $L \subset \BZ_{\ge 0} \times W$ to consists of pairs $(l, \omega)$ such that $l \ge \ell(\omega)$. Choose an arbitrary well-ordering of $W$, and equip $L$ with the lexicographic order. This is a well-ordering on $L$. 
	
	Define $r : X \to L$ as follows: $r([1, n], w, J) := \Big( \sum_{i \in [1, n]} \ell(w_i), \prod_{i \in [1, n]} w_i \Big)$. \\
	The first and second coordinates are called the \emph{length} and \emph{product} of the word, respectively. 
	
	\subsubsection{Demazure category} We will apply~\ref{rigidcat-deform}.  
	
	Define the non-full subcategory $\ms{Emb}^{\ms{d}} \subset \ms{Emb}^{\ms{c}}$ as follows: a map $([1, n], w, J) \to ([1, n], w', J)$ lies in $\ms{Emb}^{\ms{d}}$ if and only if $w_i \preceq w'_i$ for all $i$. 
	
	To define the turning map, take a solid $\ms{Dom}^{\ms{c}}$-map as shown 
	\begin{cd}[column sep = 0.9in] 
		([1,n],w,J) \ar[r, twoheadrightarrow, "{\psi : [1,n] \to [1,n']}"] & ([1,n'],w',J') \ar[r, dashed, hookrightarrow] & ([1,n'],w'',J')
	\end{cd}
	and produce the dashed $\ms{Emb}^{\ms{d}}$-map by defining $w''_j$ to be the Demazure product of the subsequence $(\psi^{-1}(j), w)$, for each $j \in [1,n']$. This lies in $\ms{Emb}^{\ms{d}}$ because 
	\[
		w'_j = \big(\text{Coxeter product of } (\psi^{-1}(j), w)\big) \preceq \big(\text{Demazure product of } (\psi^{-1}(j), w)\big) = w''_j.
	\]
	
	The first axiom expresses the associativity of the Demazure product. The second axiom expresses the Bruhat-monotonicity of the Demazure product.\footnote{This means that, if two sequences $([1, n], w)$ and $([1, n], w')$ satisfy $w_i \preceq w'_i$ for all $i$, then the Demazure product of $w$ is $\preceq$ the Demazure product of $w'$.} 
	
	Then $\ms{Rig}^{\ms{d}}$ is called the \emph{Demazure category}, because it encodes the Demazure product. 
	
	\begin{lem} \label{good}
		The function $r : X \to L$ is good in the sense of~\ref{rigidcat-deform}. 
	\end{lem}
	\begin{proof}
		First, $r$ is strictly increasing along non-identity morphisms in $\ms{Emb}^{\ms{d}}$ because, if $w_i \prec w'_i$, then $\ell(w_i) < \ell(w'_i)$. Next, fix a $\ms{Dom}^{\ms{d}}$-map as above. For each $j \in [1, n']$, basic properties of the Coxeter and Demazure products imply that 
		\[
			\ell(w'_j) \le \ell(w''_j) \le \sum_{i \in \psi^{-1}(j)} \ell(w_j), 
		\]
		and the following are equivalent: 
		\begin{itemize}
			\item At least one equality holds. 
			\item Both equalities hold. 
			\item $w'_j = w''_j$. 
			\item The sequence $(\psi^{-1}(j), w)$ refines to a reduced sequence of simple reflections.
			
			In this case, we say that $(\psi^{-1}(j), w)$ is \emph{reduced}. 
		\end{itemize}
		Summing over $j \in [1, n']$ yields the desired properties concerning $r$ and $\ms{Dom}^{\ms{d}}$. 
	\end{proof}
	
	\begin{rmk*}
		The proof shows that a morphism in $\ms{Dom}^{\ms{c}}$ or $\ms{Dom}^{\ms{d}}$ is level if and only if each $(\psi^{-1}(j), w)$ is reduced. 
	\end{rmk*}

	\mysubsec{Concrete interpretations}  \label{concrete} 
	
	\subsubsection{Direct construction} 
	The Coxeter category $\ms{Rig}^{\ms{c}}$ can be described as follows: 
	\begin{itemize}
		\item An object is a word.
		\item A morphism $([1, n], w, J) \to ([1, n'], w', J')$ is a weakly-increasing map $\psi : [1,n] \to [1,n']$ such that, for each $j \in [1, n']$, we have $J'_j \supseteq \bigcup_{i \in \psi^{-1}(j)} J_i$. 
	\end{itemize}
	The subcategory $\ms{Rig}^{\ms{d}} \subset \ms{Rig}^{\ms{c}}$ is characterized by the following requirement: 
	\begin{itemize}
		\item The above morphism belongs to $\ms{Rig}^{\ms{d}}$ if and only if, for each $j \in [1,n']$, we have 
		\[
			\big(\text{Demazure product of } (\psi^{-1}(j), w)\big) \preceq w'_j. 
		\]
	\end{itemize}
	The subcategory $\ms{Dom}^{\ms{d}} \subset \ms{Rig}^{\ms{d}}$ is characterized by the requirement that equality holds. 
	
	\subsubsection{Monoidal colimits} 
	Recall that $\mc{P}_{I, \mr{fin}}$ is the poset of finite type subsets of $I$. Consider the functor $W_\bullet : \mc{P}_{I, \mr{fin}} \to \ms{Monoids}$ defined by $J \mapsto W_J$. Consider the problem of computing $\colim W_\bullet$. The standard `amalgamation' formula for computing colimits in $\ms{Monoids}$ essentially produces the category $\ms{Dom}^{\ms{c}}$ and gives an isomorphism of monoids 
	\[
		\big(\text{connected components of }\ms{Dom}^{\ms{c}}\big) \simeq \colim W_\bullet. 
	\]
	The multiplication on the left hand side is given by concatenation of words. 
	
	Similarly, consider the functor $W^{\ms{d}}_\bullet : \mc{P}_{I, \mr{fin}} \to \ms{Monoids}$ defined by $J \mapsto W^{\ms{d}}_J$, where the latter is $W_J$ equipped with the Demazure product. Again, we have 
	\[
		\big(\text{connected components of }\ms{Dom}^{\ms{d}}\big) \simeq \colim W^{\ms{d}}_\bullet. 
	\]
	
	\begin{rmk*}
		It is a well-known idea that the Demazure product is a `deformation' of the Coxeter product, and the Hecke algebra gives an interpolation from the latter ($q = 1$) to the former ($q = 0$). This is why we chose the term `deformation' in~\ref{rigidcat-deform}. 
	\end{rmk*}
	
	\subsubsection{Demazure varieties}  \label{firstdem} 
	Let $G$ be a Kac--Moody group ind-scheme corresponding to $(I, W)$, and let $B \subset G$ be its standard Borel subgroup. For every finite type $J \subseteq I$, we have the standard type-$J$ parabolic subgroup $P_J \subset G$. For every element $w \in W$, we define $P_w \subset G$ to be the closure of the Bruhat cell $BwB$. If $w \in W_J$, then $P_w \subseteq P_J$. Equality holds if and only if $w$ is the longest element of $W_J$. 
	
	Define two functors as follows: 
	\e{
		D : \ms{Rig}^{\ms{c}} \to \ms{Varieties} & \\
		D^{\ms{d}} : \ms{Rig}^{\ms{d}} \to \ms{Varieties} &
	} 
	\begin{itemize}
		\item $D([1, n], w, J) := P_{J_1} \overset{B}{\times} \cdots \overset{B}{\times} P_{J_n} / B$. 
		\item $D^{\ms{d}}([1, n], w, J) := P_{w_1} \overset{B}{\times} \cdots \overset{B}{\times} P_{w_n} / B$. 
	\end{itemize}
	The behavior of these functors on morphisms is given by the embeddings $P_J \subset P_{J'}$ for $J \subset J'$ and the group multiplication maps $P_J \times P_J \to P_J$. There is a natural transformation
	\[
		D^{\ms{d}} \hra D|_{\ms{Rig}^{\ms{d}}}
	\]
	along closed embeddings of varieties. The values of $D$ are called Demazure (or Bott-Samelson-Demazure-Hansen) varieties, and the values of $D^{\ms{d}}$ are twisted products of Schubert varieties. 
	
	\begin{rmk*}
		The functorial factorization of each $\ms{Rig}^{\ms{d}}$-map into a $\ms{Dom}^{\ms{d}}$-map followed by an $\ms{Emb}^{\ms{d}}$-map corresponds to the functorial factorization of each map of schemes through the scheme-theoretic image. This is true because the Demazure product governs the image of Schubert varieties under convolution. Thus, $\ms{Rig}^{\ms{d}}$ (but not $\ms{Rig}^{\ms{c}}$) is a small part of an analogous construction for arbitrary schemes, see Example~\ref{rigidcat-gen}. 
	\end{rmk*}

	\mysubsec{Homotopical deletion} 
	
	\subsubsectiona 
	Fix $j = (l, \omega) \in L$, and define the full subcategory $\ms{Dom}^{\ms{c}}_{\le l, \omega} \subset \ms{Dom}^{\ms{c}}$ to consist of words with length $\le l$ and product $= \omega$. We will prove the following: 
	
	\begin{thm*} \label{deletionthm} 
		$\ms{Dom}^{\ms{c}}_{\le l, \omega}$ is contractible. 
	\end{thm*}
	
	\subsubsection{$\infty$-bistratified colimits}  
	Before beginning the proof, let us explain why this theorem is useful for computing colimits indexed by $\ms{Rig}^{\ms{d}}$. This uses everything we have done so far. 
	
	Given $F : \ms{Rig}^{\ms{d}} \to \mc{E}$, we get a chain of isomorphisms in $\tilde{\mc{E}}$ as shown: 
	\begin{cd}
		\displaystyle\colim_{\ms{Rig}^{\ms{d}}_{\le j}} \tilde{F} \ar[dash]{r}{\sim}[swap]{\ref{bis-downward}} & 
		\displaystyle\colim_{\ms{Rig}^{\ms{d}, \downarrow}_{\le j}} \tilde{F} \ar[dash]{r}{\sim}[swap]{\ref{rigid-rel}} &
		\displaystyle\colim_{\ms{Dom}^{\ms{d}}_{\le j}} \tilde{F} \ar[dash]{r}{\sim}[swap]{\ref{rigid-anti}} &
		\displaystyle\colim_{\ms{Dom}^{\ms{c}}_{\le j}} \tilde{F}^{\ms{c}} \ar[dash]{r}{\sim} &
		\displaystyle\colim_{\ms{Dom}^{\ms{c}}_{\le l, \omega}} \tilde{F}^{\ms{c}}
	\end{cd}
	The first three isomorphisms follow from the indicated references, and the fourth holds because $\ms{Dom}^{\ms{c}}_{\le l, \omega} \subset \ms{Dom}^{\ms{c}}_{\le j}$ is a connected component, and all other connected components are sent to $\hat0$ by $\tilde{F}$, so they do not affect the colimit. 
	
	\begin{lem*} \label{smooth} 
		For the following statements, we have (3) $\Leftrightarrow$ (2) $\Rightarrow$ (1). 
		\begin{enumerate}[label=(\arabic*)]
			\item $\tilde{F}$ sends basic level rank-$j$ morphisms to isomorphisms in $\tilde{\mc{E}}$. 
			\item For each basic level rank-$j$ morphism $x_1 \to x_2$, the following is a pushout square in ${\mc{E}}$
			\begin{cd}
				\displaystyle\colim_{\la \ms{Rig}^{\ms{d}}_{<j} \xra{\ms{Rig}^{\ms{d}}} x_1 \ra} F \circ \ms{tail}  \ar[r] \ar[d] & F(x_1) \ar[d] \\
				\displaystyle\colim_{\la \ms{Rig}^{\ms{d}}_{<j} \xra{\ms{Rig}^{\ms{d}}} x_2 \ra} F \circ \ms{tail}  \ar[r] & F(x_2) 
			\end{cd}
			\item For each basic level rank-$j$ morphism $x_1 \to x_2$, the following is a pushout square in ${\mc{E}}$
			\begin{cd}
				\displaystyle\colim_{\la \ms{Emb}^{\ms{d}}_{<j} \xra{\ms{Emb}^{\ms{d}}} x_1 \ra} F \circ \ms{tail}  \ar[r] \ar[d] & F(x_1) \ar[d] \\
				\displaystyle\colim_{\la \ms{Emb}^{\ms{d}}_{<j} \xra{\ms{Emb}^{\ms{d}}} x_2 \ra} F \circ \ms{tail}  \ar[r] & F(x_2) 
			\end{cd}
		\end{enumerate}
	\end{lem*}
	\begin{proof}
		The statement (2) $\Rightarrow$ (1) follows from the definition of $\tilde{F}$ as a pushout, see~\ref{bis-relative}. To show (3) $\Leftrightarrow$ (2), note that the embedding 
		\[
			\la \ms{Emb}^{\ms{d}}_{<j} \xra{\ms{Emb}^{\ms{d}}} x_1 \ra \subset \la \ms{Rig}^{\ms{d}}_{<j} \xra{\ms{Rig}^{\ms{d}}} x_1 \ra
		\]
		is a right adjoint, with left adjoint given by the functorial factorization of each $\ms{Rig}^{\ms{d}}$-map into a $\ms{Dom}^{\ms{d}}$-map followed by an $\ms{Emb}^{\ms{d}}$-map, so the embedding is final. 
	\end{proof} 
	
	Let us assume (1). Combining the theorem with the reasoning of~\ref{cases} implies 
	\begin{enumerate}[label=(\roman*)]
		\item If $l = \ell(\omega)$, then we have a pushout square 
		\begin{cd}
			\displaystyle\colim_{\la \ms{Emb}^{\ms{d}}_{<j} \xra{\ms{Emb}^{\ms{d}}} x \ra} F \circ \ms{tail}  \ar[r] \ar[d] & F(x) \ar[d] \\
			\displaystyle\colim_{\ms{Rig}^{\ms{d}}_{< j}} F \ar[r] & \displaystyle\colim_{\ms{Rig}^{\ms{d}}_{\le j}} F
		\end{cd}
		for any rank-$j$ object $x$. 
		\item If $l > \ell(\omega)$, then $\colim_{\ms{Rig}^{\ms{d}}_{< j}} F \simeq \colim_{\ms{Rig}^{\ms{d}}_{\le j}} F$. 
	\end{enumerate}
	As $j$ varies, this allows us to inductively compute $\colim_{\ms{Rig}^{\ms{d}}} F$. 
	
	The rest of this subsection is devoted to proving the theorem. 
	
	\subsubsection{Reformulation via galleries} 
	We will recast~$\ms{Dom}^{\ms{c}}_{\le l, \omega}$ in geometric terms. 
	
	The Coxeter system $(I, W)$ determines a hyperplane arrangement in a vector space $\mf{h}$, which is locally finite in the Tits cone $T \subseteq \mf{h}$. The hyperplanes in the arrangement are called \emph{walls}. Each wall partitions $T$ into three parts, obtained by intersecting with the two open half-spaces and the wall itself. Intersecting these partitions for all walls gives the partition of $T$ into \emph{faces}. A \emph{chamber} is a maximal-dimensional face. Let $\ms{Faces}$ (resp.\ $\ms{Chambers}$) denote the set of faces (resp.\ chambers). Note that $\ms{Faces}$ is a poset, where $F_1 \to F_2$ means $F_1 \subseteq \ol{F}_2$. 
	
	\begin{rmk*}
		We emphasize that faces are defined only inside $T$, and we work only inside $T$. 
	\end{rmk*}
	
	Fix a fundamental chamber $C_0$. The map $w \mapsto w\cdot C_0$ gives a bijection $W \simeq \ms{Chambers}$. Each finite type subset $J \subseteq I$ corresponds to a type-$J$ face $F_{0, J}$ of $C_0$, which has dimension $|I| - |J|$. The map $w \mapsto w \cdot F_{0, J}$ gives a bijection from $W / W_J$ to the set of type-$J$ faces. 
	
	Let $\ms{Gal}_{\le l, \omega}$ be the category defined as follows. 
	\begin{itemize}
		\item An object is a triple $([n], c, f)$, where $c = (c_0, \ldots, c_n)$ is a sequence of chambers, and $f = (f_1, \ldots, f_n)$ is a sequence of faces. We require that 
		\begin{itemize}
			\item For each $i \in [1, n]$, we have $f_i \subseteq \ol{c}_{i-1} \cap \ol{c}_i$. 
			\item $c_0 = C_0$ and $c_n = \omega \cdot C_0$. 
			\item Define the \emph{length} of the sequence $c$ to be the number of walls which it crosses, counted with multiplicity. We require that the length is $\le l$. 
		\end{itemize}
		\item A morphism $([n], c, f) \to ([n'], c', f')$ is a map $\varphi : [n'] \to [n]$ such that 
		\begin{itemize}
			\item $\varphi$ is weakly-increasing and bound-preserving, i.e.\ $\varphi(0) = 0$ and $\varphi(n') = n$. 
			\item For all $j \in [n']$, we have $c'_j = c_{\varphi(j)}$. 
			\item For all $j \in [1, n']$ and $i \in [\varphi(j-1)+1, \varphi(j)]$, we have $f'_j \subseteq \ol{f}_i$. 
		\end{itemize}
	\end{itemize}
	The objects could be called \emph{jointed generalized galleries}, but we will call them `galleries' as an abbreviation. 
	
	\begin{lem*}
		There is an isomorphism of categories $\ms{Dom}^{\ms{c}}_{\le l, \omega} \simeq \ms{Gal}_{\le l, \omega}$. 
	\end{lem*}
	\begin{proof}
		For objects, the isomorphism sends $([1, n], w, J) \mapsto ([n], c, f)$ where 
		\e{
			c_i &:= (w_1 \cdots w_i) \cdot C_0 \\
			f_i &:= (w_1 \cdots w_i) \cdot F_{0, J_i} = (w_1 \cdots w_{i-1}) \cdot F_{0, J_i}. 
		}
		The equality follows from $w_i \in W_{J_i}$. The isomorphism sends the $\ms{Dom}^{\ms{c}}_{\le l, \omega}$-morphism
		\[
			([1, n], w, J) \xra{\psi : [1,n] \to [1,n']} ([1,n'], w', J')
		\]
		to the $\ms{Gal}_{\le l, \omega}$-morphism $([n], c, f) \xra{\varphi : [n'] \to [n]} ([n'], c', f')$ where $\varphi(j)$ is defined to be the maximal index $i \in [1, n]$ such that $\psi(i) \le j$, if such an index exists, and $\varphi(j) := 0$ otherwise. 
		
		For objects, the inverse isomorphism sends $([n], c, f) \mapsto ([1, n], w, J)$, where $w_i$ is determined by $w_ic_{i-1} = c_i$, and $J_i$ is the type of the face $f_i$. The inverse isomorphism sends the $\ms{Gal}_{\le l, \omega}$-morphism 
		\[
			([n], c, f) \xra{\varphi : [n'] \to [n]} ([n'], c', f')
		\]
		to the $\ms{Dom}^{\ms{c}}_{\le l, \omega}$-morphism $([1, n], w, J) \xra{\psi : [1,n] \to [1,n']} ([1,n'], w', J')$ where $\psi(i)$ is the unique index $j$ such that $i \in [\varphi(j-1)+1, \varphi(j)]$. 
	\end{proof}
	
	\begin{rmk*}
		Here is a more conceptual description of the relationship between $\varphi$ and $\psi$. Any weakly-increasing bound-preserving map of ordered sets $\varphi : [n'] \to [n]$ induces a `pullback' map $\psi : \ms{Gaps}([n]) \to \ms{Gaps}([n'])$, where a `gap' is a pair of consecutive elements. There is an identification $\ms{Gaps}([n]) \simeq [1, n]$, where the gap $i-1 < i$ corresponds to $i \in [1, n]$. 
	\end{rmk*}
	
	\subsubsectiona 
	Define the full subcategory $\ms{Gal}^{\ms{sf}}_{\le l, \omega} \subset \ms{Gal}_{\le l, \omega}$ by requiring that each $f_i$ is not a chamber. 
	
	\begin{lem*}
		This embedding is a homotopy equivalence. 
	\end{lem*}
	\begin{proof}
		In fact, the embedding is a left adjoint, whose right adjoint sends $([n], c, f) \in \ms{Gal}_{\le l, \omega}$ to the gallery in $\ms{Gal}^{\ms{sf}}_{\le l, \omega}$ obtained by deleting the entries $(f_i, c_i)$ whenever $f_i = c_i$. 
	\end{proof}
	
	\subsubsection{From transverse paths to galleries}  \label{pg} 
	We abbreviate `piecewise linear' as PL. A PL path $\gamma : [0, 1] \to T$ is \emph{transverse} if there are only finitely many $t \in [0, 1]$ such that $\gamma(t)$ lies in a wall. Given a transverse path $\gamma$, let 
	\[
		(c_0, f_1, c_1, f_2, c_2, \ldots, f_n, c_n)
	\]
	be the sequence of faces which it passes through. As the notation suggests, this sequence of faces determines a gallery $([n], c, f)$ such that each $f_i$ is not a chamber. Define $\ms{gal}(\gamma) := ([n], c, f)$. Fix points $p_0 \in C_0$ and $p_1 \in \omega \cdot C_0$. We say that $\gamma$ is \emph{good} if its endpoints are $p_0$ and $p_1$, and $\ms{gal}(\gamma)$ has length $\le l$. In this case, we have $\ms{gal}(\gamma) \in \ms{Gal}^{\ms{sf}}_{\le l, \omega}$. 
	
	We also want to work with families of transverse paths. Thus, fix a simplicial set $K$, let $|K|$ be its geometric realization, and choose a PL map 
	\[
		\eta : |K| \times [0, 1] \to T. 
	\]
	For each subset $D \subseteq |K|$, define $\eta_D$ by restricting $\eta$ to $D \times [0, 1]$. We say that $\eta$ is a \emph{family of good transverse paths} if $\eta_{\{x\}}$ is transverse and good for each $x \in |K|$. Assume that this is the case. We will produce a map from a `refinement' of $K$ to $\ms{Gal}^{\ms{sf}}_{\le l, \omega}$. 
	
	For any subset $D \subseteq |K|$, we say that $\eta$ is \emph{constant} on $D$ if there exists a PL isomorphism $D \times [0, 1] \xra{\alpha} D \times [0, 1]$, which respects the projections to $D$ and preserves $0$ and $1$, and a transverse path $\gamma : [0, 1] \to T$, such that the following diagram commutes: 
	\begin{cd}
		D \times [0, 1] \ar[rd, swap, "\eta_D"] \ar[rr, "\alpha"] & & D \times [0, 1] \ar[ld, "\gamma \circ \pr_2"] \\
		& T
	\end{cd}
	In this case, define $\ms{gal}_\eta(D) := \ms{gal}(\gamma)$. This does not depend on the choice of $\alpha$ and $\gamma$. Let us denote this gallery by $([n], c, f)$ as above. For each $i \in [n]$, define $\eta_D^{\ms{pre}}(c_i) \subset D \times [0, 1]$ to be the set of points which are associated to the chamber $c_i$. More precisely, for any choice of $\alpha$ and $\gamma$, the index $i$ determines an open subinterval $I_i \subset [0, 1]$ such that $\gamma(I_i) \subset c_i$, and we define $\eta_D^{\ms{pre}}(c_i) := \alpha^{-1}(D \times I_i)$.\footnote{Note that $\eta_D^{\ms{pre}}(c_i) \subseteq \eta_D^{-1}(c_i)$ and the containment is strict if and only if $c_i = c_{i'}$ for some $i' \neq i$.}
	
	The nondegenerate simplices of $K$ induce a \emph{triangulation} of $|K|$, i.e.\ a partition into cells. Choose a finer triangulation of $|K|$ such that $\eta$ is constant on each cell. Let $\ms{Cells}(|K|)$ be the poset of cells, ordered by reverse closure. We have already defined the map 
	\[
		\ms{gal}_\eta : \ms{Cells}(|K|) \to \ms{Gal}^{\ms{sf}}_{\le l, \omega}
	\]
	at the level of objects, and we now enhance it to a functor. Suppose we are given a morphism in $\ms{Cells}(|K|)$, i.e.\ a containment $\ol{D} \supseteq D'$. We must define a map 
	\[
		([n], c, f) := \ms{gal}_\eta(D) \to \ms{gal}_\eta(D') =: ([n'], c', f'), 
	\]
	which is specified by a map $\varphi : [n'] \to [n]$. The definition of transverse paths implies that, for each $j \in [n']$, there exists a unique $i \in [n]$ such that $\eta_{D'}^{\ms{pre}}(c'_j) \subseteq \ol{\eta_D^{\ms{pre}}(c_i)}$. Define $\varphi(j) := i$. It is easy to check that this defines a valid morphism in $\ms{Gal}^{\ms{sf}}_{\le l, \omega}$, and that $\ms{gal}_\eta$ becomes a functor. The key idea is that each subset $\eta_{D}^{\ms{pre}}(c_i) \subseteq |K| \times [0, 1]$ admits an open neighborhood which does not intersect the preimage of any wall. 
	
	\subsubsection{From galleries to transverse paths}  \label{gp} 
	Conversely, suppose we are given a simplicial set $K$ and a map $\phi : K \to \ms{Gal}^{\ms{sf}}_{\le l, \omega}$. We will produce a family of good transverse paths $\eta$ as above. 
	
	First, note that a gallery $([n], c, f)$ can be re-encoded as a pair $([n], \mb{f})$ as follows: 
	\begin{itemize}
		\item Define an interval $I \subseteq [n]$ to be \emph{small} if it has size 1 or 2. Let $\ms{sInt}([n])$ be the poset of small intervals. 
		\item We require that $\mb{f} : \ms{sInt}([n]) \to \ms{Faces}^{\ms{op}}$ is a map of posets which sends every size-1 interval to a chamber. 
		\item Of course, we continue to require that the gallery length is $\le l$ and the start and end chambers are $C_0$ and $\omega \cdot C_0$, respectively. 
	\end{itemize}
	The values of $\mb{f}$ on size-1 (resp.\ size-2) intervals are given by $c$ (resp.\ $f$). 
	
	For each vertex $x \in K_0$, write $\phi(x) = ([n^x], \mb{f}^x)$. For each $I \in \ms{sInt}([n^x])$, choose points $t^x_I \in [0, 1]$ and $p^x_I \in \mb{f}^x_I$ ensuring that 
	\begin{itemize}
		\item $0 = t^x_{\{0\}} < t^x_{[0, 1]} < t^x_{\{1\}} < \cdots < t^x_{\{n\}} = 1$. 
		\item $p^x_{\{0\}} = p_0$ and $p^x_{\{n\}} = p_1$. 
	\end{itemize}
	For constructing $\eta$, our main idea is to require that $\eta_x$ sends $t^x_I \mapsto p^x_I$ and then to extend $\eta_x$ to $[0, 1]$ by linear interpolation. Clearly, this path is transverse and good. 
	
	Next, we want to define $\eta$ on all of $|K| \times [0, 1]$ using linear interpolation. To do so, we first define a triangulation of $|K| \times [0, 1]$ whose 0-cells (i.e.\ vertices) are given by the points $(x, t^x_I)$ for $x \in K_0$ as above. We will do this over each simplex $\sigma : \Delta^m \to K$ and then patch the simplices together. The simplex $\sigma$ gives a collection of galleries and maps: 
	\e{
		\phi(\sigma(i)) &= ([n^{\sigma(i)}], \mb{f}^{\sigma(i)}) \\
		\phi(\sigma(i \to j)) &= \Big( \phi(\sigma(i)) \xra{\varphi^{j, i} : [n^{\sigma(j)}] \to [n^{\sigma(i)}]} \phi(\sigma(j)) \Big). 
	}
	Here $i \le j$ are elements of $[m]$. Define a poset $\mc{P}_\sigma$ as follows: 
	\begin{itemize}
		\item An object is a pair $(i, I)$ where $i \in [m]$ and $I \in \ms{sInt}([n^{\sigma(i)}])$. 
		\item There is a map $(i, I) \to (j, J)$ if and only if $i \le j$ and $I \subseteq \varphi^{j, i}_*(J)$. 
		
		By definition, $\varphi^{j, i}_*(J)$ is the interval in $[n^{\sigma(i)}]$ defined by the $\varphi^{j, i}$-images of the endpoints of $J$. It need not be small. 
	\end{itemize}
	It is easy to show that the geometric realization (of the nerve) of this poset, denoted $|\mc{P}_\sigma|$, is homeomorphic to $\Delta^m \times [0, 1]$, via the map defined by $(i, I) \mapsto (i, t^{\sigma(i)}_I)$ and linear interpolation. (The most efficient way to prove this is to use induction on $m$. In the base case $m = 0$, the poset $\mc{P}_\sigma = \ms{sInt}([n^{\sigma(0)}])$ is a zig-zag, so its geometric realization is homeomorphic to $[0, 1]$.) Now we may define $\eta$ over $\Delta^m$ by linear interpolation within the simplices of $\mc{P}_\sigma$. This construction is functorial with respect to $\sigma$, so it defines $\eta$ on all of $|K| \times [0, 1]$. 
	
	It is straightforward to check that $\eta$ is a family of good transverse paths. Furthermore, $\eta$ is constant on the interior of each simplex of $K$. Thus, the already-established `paths to galleries' direction turns $\eta$ into a functor $\ms{gal}_\eta : \ms{Cells}^{\ms{nondeg}}(|K|) \to \ms{Gal}^{\ms{st}}_{\le l, \omega}$ whose domain is the set of nondegenerate simplices of $K$. By construction, there is a commutative diagram 
	\begin{cd}
		K \ar[rr, leftarrow, "\ms{ev}_0"] \ar[rd, swap, "\phi"] & & \ms{Cells}^{\ms{nondeg}}(|K|) \ar[ld, "\ms{gal}_\eta"] \\
		& \ms{Gal}^{\ms{sf}}_{\le l, \omega} 
	\end{cd}
	where the map $\ms{ev}_0$ sends a nondegenerate simplex $\sigma : \Delta^m \to K$ to the vertex $\sigma(0) \in K$. 
	
	\subsubsection{Prefix-straightening} 
	As before, suppose we are given a simplicial set $K$ and a map $\phi : K \to \ms{Gal}^{\ms{sf}}_{\le l, \omega}$. We will show that, after modifying $K$ via some homotopy equivalences, the map $\phi$ becomes constant. This implies that $\ms{Gal}^{\ms{sf}}_{\le l, \omega}$ is contractible. 
	
	Use the `galleries to paths' direction to construct a family $\eta : |K| \times [0, 1] \to T$. Extend this to a family ${h} : \big(|K| \times [0, 1]\big) \times [0, 1] \to T$ as follows: 
	\[
		{h}(x, s, t) := \begin{cases}
			\frac{t}{s} \cdot p_0 + (1 - \frac{t}{s}) \cdot \eta(x, s) & \text{ if } t \le s \\
			\eta(x, t) & \text{ otherwise} 
		\end{cases}
	\]
	Intuitively, at $s \in [0, 1]$, each path $\eta_x : [0, 1] \to T$ is modified by replacing the $[0, s]$-subpath by a straight line segment with the same endpoints. This is why we chose the term `prefix-straightening.' At $s = 0$, these paths agree with $\eta$, but at $s = 1$ these paths are all equal to the straight line segment from $p_0$ to $p_1$. These modified paths are transverse because the straight line segment on $[0, s]$ contains the point $p_0$ which does not lie in any wall. 
	
	Next, the `paths to galleries' direction gives the following maps: 
	\begin{cd}[column sep = 6mm]
		K \ar[r, leftarrow, "\ms{ev}_0"] \ar[rrd, swap, "\phi"] & \ms{Cells}^{\ms{nondeg}}(|K|) \ar[r, leftarrow] \ar[rd, "\ms{gal}_\eta"] & \ms{Cells}(|K| \times \{0\}) \ar[r, hookrightarrow] \ar[d] & \ms{Cells}(|K| \times [0, 1]) \ar[r, hookleftarrow] \ar[ld, swap, "\ms{gal}_{{h}}"] & \ms{Cells}(|K| \times \{1\}) \ar[lld] \\
		& & \ms{Gal}^{\ms{sf}}_{\le l, \omega} 
	\end{cd}
	Here $\ms{Cells}(|K| \times [0, 1])$ refers to an arbitrary triangulation which satisfies the following: it works for ${h}$, the subspaces $|K| \times \{0\}$ and $|K| \times \{1\}$ are unions of cells, and the induced triangulation of $|K| \times \{0\}$ refines the `nondegenerate simplex' triangulation of $|K|$. Then the second horizontal map, which sends each such cell of $|K| \times \{0\}$ to the nondegenerate simplex which contains it, is well-defined. 
	
	After possibly modifying $K$ via a homotopy equivalence, we may assume that every nondegenerate simplex $\Delta^m \to K$ is a monomorphism. This implies that $\ms{ev}_0$ is a homotopy equivalence. Clearly, the other horizontal arrows are homotopy equivalences as well. The last downward map is constant, with value equal to the gallery associated to the path given by the straight line segment from $p_0$ to $p_1$. This concludes the proof of Theorem~\ref{deletionthm}. 
	
	\begin{rmk*}
		We can now explain the term `homotopical deletion.' There is a standard algorithm for turning a nonreduced sequence of simple reflections in $W$ into a reduced sequence with the same product. Namely, take the smallest initial subsequence which is nonreduced, modify it using the Deletion Property~\cite[Prop.\ 1.4.7]{comb-cox} to make it reduced, and continue until the whole sequence is reduced. The `prefix-straightening' idea is a geometric analogue of this algorithm, which has the advantage of working nicely in families, i.e.\ it is `homotopical.'
	\end{rmk*}

	\section{Colimit presentation of the Hecke \texorpdfstring{$\infty$}{oo}-category}  \label{s-main} 
	
	For psychological ease, we work with an algebraic loop group, even though all of the subsequent material works for any simply-connected Kac--Moody group. 
	
	Let $G$ be a simply-connected semisimple algebraic group over a field of characteristic zero, let $\LG$ be its loop group, and let $I$ be the Coxeter--Dynkin diagram of $\LG$. Then the affine Coxeter group $(I, W)$ is the Weyl group of $\LG$. Let $\mb{I} \subset \LG$ denote the standard Iwahori subgroup, and let $\mb{P}_J \subset \LG$ be the standard parahoric subgroup of type $J \subset I$. (Thus $\mb{P}_I = \LG$ and $\mb{P}_\emptyset = \mb{I}$.) 
	
	Let $\mc{P}_{I, \mr{fin}}$ be the poset of finite type subsets of $I$. In~\ref{conv}, we construct a functor 
	\[
		\mc{H}_\bullet : \mc{P}_{I, \mr{fin}} \to \alg(\cat)
	\]
	which sends $J \subset I$ to the Hecke category $\Shv(\mb{I} \backslash \mb{P}_J / \mb{I})$. Our goal is to show that $\colim \mc{H}_\bullet \simeq \mc{H}_I$. In~\ref{ss-amalg}, we express the underlying DG-category of $\colim  \mc{H}_\bullet$ as $\colim \wh{\mc{H}}_\bullet$ for a functor 
	\[
		\wh{\mc{H}}_\bullet : (\mc{P}_{I, \mr{fin}})^\amalg_{\mr{act}} \to \cat
	\]
	from a larger category $(\mc{P}_{I, \mr{fin}})^\amalg_{\mr{act}}$. In~\ref{ss-interpret}, we recast this as the colimit over a functor 
	\[
		\wt{\mc{H}}_\bullet : \ms{Rig}^{\ms{d}} \to \cat
	\]
	defined on the Demazure category. We also interpret this functor geometrically, in terms of twisted products of Schubert varieties. Finally, in~\ref{ss-checkbis}, we use our earlier work on $\infty$-bistratified colimits (\ref{smooth}) to show that $\colim \wt{\mc{H}}_\bullet \simeq \mc{H}_I$, which is our main theorem.

	\mysubsec{Definition of the Hecke \texorpdfstring{$\infty$}{oo}-category} \label{conv} 
	
	\subsubsection{Overview} The $\infty$-category $\Shv(\mb{I} \backslash \LG / \mb{I})$ should have a `convolution' monoidal structure for which the tensor product is given by convolution (i.e.\ pull-push along a span) and the associativity constraint for a 3-term multiplication is given by proper base-change. In order to make this a monoidal structure in the $\infty$-categorical sense, one must provide compatibilities between these base-change isomorphisms for each multiplication involving $\ge 4$ terms. 
	
	One approach to doing so is to realize $\Shv(-)$ as a monoidal functor out of the $\infty$-category of spans. For spans of laft\footnote{Roughly speaking, a prestack is `laft' if it can be recovered from its values on finite type affine schemes.} prestacks, this has been done in~\cite[Vol.\ 1, Chap.\ 9]{gr} and applied to convolution monoidal structures on $\D$-module $\infty$-categories in~\cite[Vol.\ 2, Chap.\ 4]{gr}. Unfortunately, $\mathbf{I} \backslash \LG / \mathbf{I}$ is not laft.
	
	One possible fix is to extend the framework of \cite{gr} to include quotients of laft prestacks by pro-algebraic groups. Instead, we will use the theory of $\Shv^*(\LG)$-actions on $\infty$-categories to obtain the desired monoidal structure. What these two approaches have in common is that both involve approximating $\LG$ by schemes of finite type.
	
	We will use the paper~\cite{act} which defines two versions of $\D$-module categories on ind-pro-schemes. Another development of $\D$-modules on ind-pro-schemes is given in~\cite{raskin}. 
	
	\subsubsection{Main definition} \label{conv-def} 
	Denote the affine flag variety by $\Fl := \LG / \mb{I}$. We interpret the notation $\Shv(\mathbf{I} \backslash \LG / \mathbf{I})$ to mean the $\infty$-category ${}_{\mathbf{I}}\Shv^*(\Fl)$, where the $\mathbf{I}$-subscript means coinvariants in the sense of~\cite[4.1.1]{act}. 
	
	\subsubsectiona \label{conv1} 
	To see that this definition is reasonable, let us describe ${}_{\mathbf{I}}\Shv^*(\Fl)$ more concretely. Let $\Fl \simeq \colim_n \Fl^{\le n}$ be an ind-scheme presentation of $\Fl$ such that each stratum $\Fl^{\le n}$ is $\mathbf{I}$-invariant. Then
	\e{
		{}_{\mathbf{I}}\Shv^*(\Fl) &\simeq {}_{\mathbf{I}}\big(\! \colim_n \Shv(\Fl^{\le n}) \big) \\
		&\simeq \colim_n {}_{\mathbf{I}}\Shv(\Fl^{\le n})
	}
	where the first line follows by definition of $\Shv^*(-)$ on ind-schemes, and the second line follows because coinvariants are defined by colimits. The functoriality in the colimit diagram is via $*$-pushforward. The $\infty$-categories ${}_{\mathbf{I}}\Shv(\Fl^{\le n})$ can be understood as follows: 
	
	\begin{lem*} 
		For a fixed $n$, let $\mathbf{I}' \hra \mathbf{I}$ be a normal pro-unipotent sub-pro-group such that
		\begin{itemize}
			\item The action of $\mathbf{I}'$ on $\Fl^{\le n}$ is trivial.
			\item The quotient $Q := \mathbf{I} / \mathbf{I}'$ is finite type.
		\end{itemize}
		Then ${}_{\mathbf{I}}\Shv(\Fl^{\le n}) \simeq \Shv(Q \backslash \Fl^{\le n})$.
	\end{lem*}
	\begin{proof}
		By definition, ${}_{\mathbf{I}}\Shv^*(\Fl^{\le n})$ is the colimit of the diagram of categories given by applying $\Shv^*(-)$ to the cosimplicial diagram for the action of $\mathbf{I}$ on $\Fl^{\le n}$. By \cite[Cor.\ 5.5.3.4]{htt} or \cite[3.1.2]{act}, a colimit of categories along left adjoints is equivalent to a limit along the corresponding right adjoints, so this is equivalent to the limit of the diagram of categories given by applying $\Shv^*(-)$ with the $*$-pullback functoriality, which exists because the action maps are smooth, see~\cite[3.1.3]{act}. Because $\mathbf{I}'$ is pro-unipotent, it is contractible, so in the action diagram we can replace $\mathbf{I}$ by $Q$. Now our diagram involves finite type schemes, so smooth descent for $\D$-modules implies that the resulting limit category is $\Shv(Q \backslash \Fl^{\le n})$.
	\end{proof}
	
	\begin{lem}\label{inv-coinv2}
		We have $\Shv^*(\LG)_{\mathbf{I}} \simeq \Shv^*(\Fl)$ as left $\Shv^*(\mathbf{I})$-modules.
	\end{lem}
	\begin{proof}
		Let $q : \LG \to \Fl$ be the quotient map. Then $\LG \simeq \colim_n q^{-1}(\Fl^{\le n})$ is an ind-pro-scheme presentation of $\LG$ such that each stratum $q^{-1}(\Fl^{\le n})$ is invariant under the right $\mathbf{I}$-action. We have
		\e{
			\Shv^*(\LG)_{\mathbf{I}} &\simeq \big( \colim_n \Shv^*(q^{-1}(\Fl^{\le n})) \big)_{\mathbf{I}} \\
			&\simeq \colim_n \Shv^*(q^{-1}(\Fl^{\le n}))_{\mathbf{I}}.
		}
		The first line follows from the definition of $\Shv^*(-)$ for ind-pro-schemes given in~\cite[3.3.3]{act}, and the second line follows from the fact that coinvariants are given by a colimit.
		
		To finish, we will prove that
		\[
		\Shv^*(q^{-1}(\Fl^{\le n}))_{\mathbf{I}} \simeq \Shv(\Fl^{\le n}).
		\]
		The left hand side is defined as the colimit of the diagram of categories given by applying $\Shv^*(-)$ to the cosimplicial diagram for the (right) action of $\mathbf{I}$ on $q^{-1}(\Fl^{\le n})$. As in Lemma~\ref{conv1}, we can replace this by the limit over the diagram obtained by applying $\Shv^*(-)$ with the $*$-pullback functoriality. Let $\mathbf{I}_r \hra \mathbf{I}$ be a sequence of progressively smaller normal pro-unipotent sub-pro-groups, such that $\mathbf{I} / \mathbf{I}_r$ is finite type for all $r$. Then, by definition of $\Shv^*(-)$, the preceding diagram is the limit over $r$ (with respect to $*$-pushforward maps) of the result of applying $\Shv(-)$ with $*$-pullback functoriality to the action diagram for $\mathbf{I} / \mathbf{I}_r$ acting on $q^{-1}(\Fl^{\le n}) / \mathbf{I}_r$. By smooth descent for $\D$-modules, the limit of the $r$-th such diagram is $\Shv(\Fl^{\le n})$ because $q^{-1}(\Fl^{\le n}) / \mathbf{I} \simeq \Fl^{\le n}$. Since limits commute with limits, we may now take the limit over $r$, and the claim follows.
	\end{proof}
	
	\subsubsectiona 
		We will construct the \emph{convolution monoidal structure} on ${}_{\mathbf{I}}\Shv^*(\Fl)$. By~\cite[Lem.\ 3.4.5]{act}, $\Shv^*(\LG)$ has a monoidal structure given by $*$-pushforward along the multiplication map. Accordingly, we may consider the $(\infty, 2)$-category of $\Shv^*(\LG)$-modules in $\cat$, denoted $\Shv^*(\LG)\mod$. By~\cite[2.3.9]{act}, we have
		\[
			\on{Hom}_{\Shv^*(\LG)\mod}(\Shv^*(\LG)_{\mathbf{I}}, \Shv^*(\LG)_{\mathbf{I}}) \simeq {}^{\mathbf{I}}\Shv^*(\LG)_{\mathbf{I}}
		\]
		This is proved by considering the cosimplicial diagram that defines $\mathbf{I}$-coinvariants. Next, since $\mathbf{I}$ is an extension of a pro-unipotent group scheme by a finite type group scheme, \cite[Thm.\ 4.2.4]{act} says that $\mathbf{I}$-invariants are equivalent to $\mathbf{I}$-coinvariants, so we have  
		\[
			{}^{\mathbf{I}}\Shv^*(\LG)_{\mathbf{I}} \simeq {}_{\mathbf{I}}\Shv^*(\LG)_{\mathbf{I}}
		\]
		Finally, Lemma~\ref{inv-coinv2} implies that
		\[
			{}_{\mathbf{I}}\Shv^*(\LG)_{\mathbf{I}} \simeq {}_{\mathbf{I}}\Shv^*(\Fl)
		\]
		The convolution monoidal structure on ${}_{\mathbf{I}}\Shv^*(\Fl)$ is the one which corresponds to the composition monoidal structure on $\on{Hom}_{\Shv^*(\LG)\mod}(\Shv^*(\LG)_{\mathbf{I}}, \Shv^*(\LG)_{\mathbf{I}})$. 
	
	\subsubsection{Finite type subcategories} \label{conv-sub} 
		For any finite type $J \subset I$, we interpret the notation $\Shv(\mb{I} \backslash \mb{P}_J / \mb{I})$ to mean the full subcategory 
		\[
			{}_{\mb{I}}\Shv(\Fl_J) \hra {}_{\mb{I}}\Shv^*(\Fl) 
		\]
		where $\Fl_J := \mb{P}_J / \mb{I}$ is a finite dimensional flag variety which embeds into $\Fl$. This is compatible with the main definition~\ref{conv-def}.
		
		To define the convolution monoidal structure on ${}_{\mb{I}}\Shv(\Fl_J)$, use \cite[Prop.\ 2.2.1.1]{ha}, which states that, if the set of objects of a full subcategory of a monoidal $\infty$-category (such as ${}_{\mb{I}}\Shv^*(\Fl)$) is invariant under the tensor product, then the full subcategory acquires a unique monoidal structure compatible with the ambient one. 
	
	\begin{rmk*}
		This definition is reasonable because ${}_\mathbf{I}\Shv(\Fl_J)$ agrees with the finite type Hecke category associated to $\mb{P}_J$. For example, suppose $\mb{P}_J \subseteq G(\oh)$, so that the first congruence subgroup $\mathbf{G}_1$ is normal in $\mb{P}_J$. Let $P_J := \mb{P}_J / \mathbf{G}_1$ and $B := \mathbf{I} / \mathbf{G}_1$, so that $\Fl_J \simeq P_J/B$. The proof of Lemma~\ref{conv1} implies that ${}_\mathbf{I}\Shv(\Fl_J) \simeq {}_{B}\Shv(\Fl_J)$ because $\mathbf{G}_1$ is pro-unipotent and acts trivially. Therefore, ${}_{\mathbf{I}}\Shv(\Fl_J) \simeq \Shv(B \backslash P_J / B)$. It is also easy to identify the monoidal structures on both sides.
	\end{rmk*}
	
	\subsubsection{Hecke functor}  \label{last} 
	We can now define the functor 
	\[
		\mc{H}'_\bullet : \mc{P}_{I, \mr{fin}}^{\triangleright} \to \alg(\cat)
	\]
	by requiring that $\mc{H}'_I := {}_{\mb{I}}\Shv^*(\Fl)$, $\mc{H}'_J := {}_{\mb{I}}\Shv(\Fl_J)$ for all finite type $J \subset I$, and the arrow $J \hra I$ goes to full embedding
	$
		{}_{\mb{I}}\Shv(\Fl_J) \hra {}_{\mb{I}}\Shv^*(\Fl) 
	$
	from~\ref{conv-sub}.  
	
	After these requirements are in place, the behavior of $\mc{H}'_\bullet$ on higher simplices and the compatibility with the monoidal structures both come for free. Indeed, the higher simplices come from the fact that the $\infty$-category of fully faithful functors to a given $\infty$-category (e.g.\ ${}_{\mb{I}}\Shv^*(\Fl)$) is equivalent to its poset of full subcategories. The compatibility with monoidal structures comes from the second paragraph of~\ref{conv-sub}. 
	
	Finally, let $\mc{H}_\bullet$ be the restriction of $\mc{H}'_\bullet$ to $\mc{P}_{I, \mr{fin}}$.

	\mysubsec{Monoidal colimits are amalgamated products} \label{ss-amalg}
	
	\subsubsection{Type sequences} \label{preamalg-def} 
	Let $(\mc{P}_{I, \mr{fin}})^{\amalg}_{\mr{act}}$ be the category defined as follows: 
	\begin{itemize}
		\item An object is a pair $([1, n], J)$, where $J = (J_1, \ldots, J_n)$ is a sequence in $\mc{P}_{I, \mr{fin}}$. 
		\item A morphism $\varphi : ([1,n], J) \to ([1,n'], J')$ is a weakly-increasing map $\varphi_* : [1,n] \to [1,n']$ such that, for each $j \in [1,n']$, we have $J'_j \supseteq \bigcup_{i \in \varphi_*^{-1}(j)} J_i$. 
	\end{itemize}
	As the notation suggests, this is a special case of a more general construction, in which $\mc{P}_{I, \mr{fin}}$ is replaced by an arbitrary poset. 
	
	\begin{prop} \label{preamalg-prop} 
		There exists a functor 
		\[
			\wh{\mc{H}}_\bullet : (\mc{P}_{I, \mr{fin}})^{\amalg}_{\mr{act}} \to \cat
		\]
		with the following properties: 
		\begin{itemize}
			\item For each $([1, n], J)$ as above, we have 
			\[
				\wh{\mc{H}}_{([1, n], J)} = \mc{H}_{J_1} \underset{\mc{H}_\emptyset}{\otimes} \cdots \underset{\mc{H}_\emptyset}{\otimes} \mc{H}_{J_n}
			\]
			\item $\wh{\mc{H}}_\bullet$ sends a morphism $\varphi$ (as above) to the functor described as follows. For each $j \in [1, n']$, write $\varphi_*^{-1}(j) = [a_j, b_j]$, and consider the functor 
			\[
				\mc{H}_{J_{a_j}} \underset{\mc{H}_\emptyset}{\otimes} \mc{H}_{J_{a_j+1}} \underset{\mc{H}_\emptyset}{\otimes} \cdots \underset{\mc{H}_\emptyset}{\otimes} \mc{H}_{J_{b_j}} \to \mc{H}_{J'_j}
			\]
			defined by the convolution monoidal structure. Then tensor over $j \in [1, n']$. 
			\item There is an equivalence of DG-categories 
			\[
				\colim_{(\mc{P}_{I, \mr{fin}})^{\amalg}_{\mr{act}}} \wh{\mc{H}}_\bullet \xra{\sim} \mr{oblv}(\colim_{\mc{P}_{I, \mr{fin}}} \mc{H}_\bullet). 
			\]
		\end{itemize}
	\end{prop}
	
	In the last bullet, `oblv' refers to forgetting the monoidal structure. The rest of the subsection is devoted to proving the proposition. 
	
	\subsubsectiona  \label{dual-incarnation}
	Following~\cite[4.1.2.8]{ha}, we think of the category $\Delta^\op$ in its `dual incarnation,' i.e.\ objects are ordered sets $\bb{n} := [1, n]$ (so $\bb{0} = \emptyset$), and a morphism $\bb{n_1} \to \bb{n_2}$ is given by an order-preserving map $\{-\infty\} \sqcup \bb{n_1} \sqcup \{\infty\} \to \{-\infty\} \sqcup \bb{n_2} \sqcup \{\infty\}$ which also preserves $\pm \infty$. The equivalence with $\Delta^\op$ is given by $\bb{n} \mapsto \Delta^n$.
	
	Using the `dual incarnation' is psychologically helpful because it agrees with the intuition that the morphisms in~\ref{preamalg-def} should be visualized as collisions of labeled points. On the other hand, we continue to denote this category by $\Delta^\op$ to ensure that our notations agree with~\cite{ha}. 
	
	\subsubsectiona  \label{ha-defs} Since we are interested in associative (not commutative) algebras, we use the theory of planar $\infty$-operads from \cite[4.1.3]{ha}. Here are the main concepts we need.
	\begin{enumerate}[label=(\roman*)]
		\item A \emph{planar $\infty$-operad} is an $\infty$-category $\oh^{\otimes}$ equipped with a functor $p : \oh^{\otimes} \to \assoc$ satisfying the conditions of~\cite[Def.\ 4.1.3.2]{ha}.
		\item A \emph{monoidal $\infty$-category} is an $\infty$-operad $p : \mc{C}^\otimes \to \ner(\Delta^\op)$ for which $p$ is a cocartesian fibration, see~\cite[Def.\ 2.1.2.13]{ha}. We abuse notation by also calling the fiber $\mc{C} := \mc{C}^{\otimes} \underset{\ner(\Delta^\op)}{\times} \{ \bb1\}$ a monoidal $\infty$-category.
		\item Given planar $\infty$-operads $\oh^\otimes$ and $\oh^{\prime\otimes}$, the \emph{$\infty$-category of $\oh^{\otimes}$-algebras in $\oh^{\prime\otimes}$} is the full subcategory of $\on{Fun}_{\assoc}(\oh^\otimes, \oh^{\prime\otimes})$ spanned by planar $\infty$-operad maps, see \cite[Def.\ 2.1.2.7]{ha}. This category is denoted $\alg_{\oh}(\oh')$. When $\oh^\otimes = \ner(\Delta^\op)$ and $\oh^{\prime\otimes}$ is a monoidal $\infty$-category, this is also called the  \emph{$\infty$-category of associative algebras in $\oh'$}.
		\item A \emph{marked simplicial set} is a pair $(X, M)$ where $X$ is a simplicial set and $M$ is a collection of `marked' edges of $X$ which includes all the degenerate edges. The \emph{category of planar $\infty$-preoperads} is the category of marked simplicial sets over $\assoc^{\natural}$ where the $\natural$ sign means that the marked edges are the inert ones. For details, see~\cite[Def.\ 3.1.0.1]{htt}, \cite[Def.\ 2.1.4.2]{ha}, and \cite[Def.\ 4.1.3.12]{ha}.
	\end{enumerate}
	
	\subsubsectiona \label{coprod} 
	We record a variant of \cite[Const.\ 2.4.3.1]{ha} which turns an $\infty$-category $\mc{I}$ into a planar $\infty$-operad $\mc{I}^\amalg$. 
		
	Let $\Delta^{\op, *}$ be the category defined as follows:
	\begin{itemize}
		\item The objects are pairs $(\bb{n}, i)$ where $\bb{n}$ is a totally ordered finite set and $i \in \bb{n}$.
		\item A morphism $\varphi : (\bb{n_1}, i_1) \to (\bb{n_2}, i_2)$ is a weakly increasing map $\varphi_* : \{-\infty\} \sqcup \bb{n_1} \sqcup \{\infty\} \to \{-\infty\} \sqcup \bb{n_2} \sqcup \{\infty\}$ which preserves $\pm \infty$ and satisfies $\varphi(i_1) = i_2$.
	\end{itemize}
	For any simplicial set $\mc{I}$, we define a simplicial set $\mc{I}^{\amalg}$ equipped with a map $p : \mc{I}^{\amalg} \to \assoc$ so that the following universal property is satisfied: for every map of simplicial sets $K \to \assoc$, we have a canonical bijection
	\[
	\Hom_{\assoc}(K, \mc{I}^\amalg) \simeq \Hom_{\sset}\Big(K \underset{\assoc}{\times} \ner(\Delta^{\op, *}), \mc{I}\Big).
	\]
	As in \cite[Rmk.\ 2.4.3.2]{ha}, we have a canonical isomorphism of simplicial sets $\mc{I}^{\amalg} \underset{\assoc}{\times} \{\bb{n}\} \simeq \mc{I}^{\times n}$. The proof of \cite[Prop.\ 2.4.3.3]{ha} implies that, if $\mc{I}$ is an $\infty$-category, then $p : \mc{I}^{\amalg} \to \assoc$ is an $\infty$-operad. As in \cite[Ex.\ 2.4.3.5]{ha}, the forgetful functor $\Delta^{\op, *} \to \Delta^\op$ induces a map of planar $\infty$-preoperads $\mc{I}^\flat \times \assoc^{\natural} \to \mc{I}^{\amalg, \natural}$. Over $\bb{n} \in \Delta^\op$, this map is the diagonal embedding $\mc{I} \to \mc{I}^{\times n}$.
	
	A morphism in $\mc{I}^\amalg$ is called \emph{active} if its image in $\Delta^\op$ is a morphism $\{-\infty\} \sqcup \bb{n_1} \sqcup \{\infty\} \to \{-\infty\} \sqcup \bb{n_2} \sqcup \{\infty\}$ which sends $\bb{n_1}$ into $\bb{n_2}$. Let $\mc{I}^\amalg_{\mr{act}} \subset \mc{I}^{\amalg}$ be the full subcategory spanned by the active morphisms. This generalizes~\ref{preamalg-def}. 
	
	\subsubsectiona The next theorem says that, in the $\infty$-category of planar $\infty$-operads, the object $\mc{I}^\amalg$ corepresents diagrams of algebras indexed by $\mc{I}$. 
	
	\begin{thm*} \label{wr}
		For any simplicial set $\mc{I}$ and any planar $\infty$-operad $\mc{O}^\otimes$, pullback along the map $\mc{I}^\flat \times \assoc^{\natural} \to \mc{I}^{\amalg, \natural}$ from~\ref{coprod} induces an equivalence of $\infty$-categories
		\[
		\on{Map}^\flat_{/\assoc^{\natural}}(\mc{I}^{\amalg, \natural}, \mc{O}^{\otimes, \natural}) \xrightarrow{\sim} \on{Map}^\flat_{/\assoc^{\natural}}(\mc{I}^\flat \times \assoc^{\natural}, \mc{O}^{\otimes, \natural}).
		\]
	\end{thm*}
	\begin{proof}
		The proof of \cite[Thm.\ 2.4.4.3]{ha} does not actually use that $\Shv^\otimes$ is an $\infty$-operad. Our desired statement follows if we apply this proof with the following replacements:
		\e{
			\on{N}(\mc{F}\mathrm{in}_*) &\rightsquigarrow \assoc \\
			\mc{C}^{\otimes, \natural} &\rightsquigarrow \assoc^{\natural} \\
			\Shv^{\otimes, \natural} &\rightsquigarrow \mc{I}^{\flat}
		}
		where the structure map $\mc{I}^\flat \to \assoc$ sends everything to $\bb{1}$. The $\infty$-preoperad $\mc{C}^{\otimes, \natural} \odot \Shv^{\otimes, \natural}$ becomes the planar $\infty$-preoperad $\mc{I}^{\flat} \times \assoc^{\natural}$ which maps to $\assoc^\natural$ via projection to the second factor. The $\infty$-preoperad $(\mc{C}^{\otimes} \wr \Shv^{\otimes}, M)$ becomes $\assoc \times_{\assoc} \mc{I}^{\amalg} \simeq \mc{I}^{\amalg}$ with marked edges given by inert morphisms in $\mc{I}^{\amalg}$. Although \cite[Thm.\ 2.4.4.3]{ha} only claims a weak equivalence of $\infty$-preoperads, its proof establishes an equivalence of mapping categories, i.e.\ between the $\infty$-categories given by $\on{Map}^\flat(-, \oh^{\otimes, \natural})$ rather than between their maximal Kan-fibrant subcomplexes $\on{Map}^{\sharp}(-, \oh^{\otimes, \natural})$.
	\end{proof}
	
	\subsubsectiona  \label{wr-cor} Fix an $\infty$-category $\mc{I}$, and let $\pi : \mc{I} \to \{*\}$ be the projection. Let $\pi' : \mc{I}^\amalg \to \{*\}^\amalg \simeq \assoc$ be the map induced by the functoriality of~\ref{coprod}.
	
	\begin{cor*}
		Let $\oh^\otimes$ be a planar $\infty$-operad. There is a commutative diagram of $\infty$-categories
		\begin{cd}
			\alg_{\assoc}(\oh) \ar[r, "\sim"] \ar[d, "\on{Res}_{\pi'}"] & \on{Fun}(\{*\}, \alg_{\assoc}(\oh)) \ar[d, "\on{Res}_\pi"] \\
			\alg_{\mc{I}^{\amalg}}(\oh) \ar[r, "\sim"] & \on{Fun}(\mc{I}, \alg_{\assoc}(\oh))
		\end{cd}
		where the horizontal functors are equivalences.
	\end{cor*}
	\begin{proof}
		The bottom horizontal functor is defined by
		\e{
			\alg_{\mc{I}^{\amalg}}(\oh) &\simeq \on{Map}^\flat_{/\assoc^{\natural}}(\mc{I}^{\amalg, \natural}, \mc{O}^{\otimes, \natural}) \\
			&\simeq \on{Map}^\flat_{/\assoc^{\natural}}(\mc{I}^\flat \times \assoc^{\natural}, \mc{O}^{\otimes, \natural}) \\
			&\simeq \on{Fun}(\mc{I}, \alg_{\assoc}(\oh)).
		}
		The first equivalence is the definition of $\alg(-)$. The second equivalence is Theorem~\ref{wr}. The third equivalence follows from the definition of the functor category and the adjunctions noted in~\cite[3.1.3]{htt}: for any simplicial set $K$,
		\e{
			\Hom_{\sset}(K, \on{Fun}(\mc{I}, \alg_{\assoc}(\oh))) &\simeq \Hom_{\sset}(K \times \mc{I}, \alg_{\assoc}(\oh)) \\
			&\simeq \Hom_{\sset}(K \times \mc{I}, \on{Map}^\flat_{/\assoc^\natural}(\assoc^\natural, \oh^{\otimes, \natural})) \\
			&\simeq \Hom_{/\assoc^\natural}(K^\flat \times \mc{I}^\flat \times \assoc^\natural, \oh^{\otimes, \natural}) \\
			&\simeq \Hom_{\sset}(K, \on{Map}^\flat_{/\assoc^{\natural}}(\mc{I}^\flat \times \assoc^{\natural}, \mc{O}^{\otimes, \natural})).
		}
		The commutative square results from functoriality in $\mc{I}$.
	\end{proof}
	
	\begin{lem} \label{fill}
		Let $q : \mc{C} \to \Shv$ be a cocartesian fibration of $\infty$-categories, and let $d \in \Shv$ be a terminal object.
		\begin{enumerate}[label=(\roman*)]
			\item Let $\mc{C}_d$ be the fiber of $q$ over $d$. The embedding $\iota : \mc{C}_d \to \mc{C}$ admits a left adjoint $\ell$.
			\item Suppose we are given a commutative diagram of $\infty$-categories
			\begin{cd}
				\mc{K} \ar[r, "p"] \ar[d, hookrightarrow] & \mc{C} \ar[d, "q"] \\
				\mc{K}^{\triangleright} \ar[r, "p'"] & \Shv
			\end{cd}
			Then there is an equivalence
			\[
			\Big(\mc{C}_{p/} \underset{\Shv_{qp/}}{\times} \Shv_{p'\!/}\Big) \underset{\mc{C}}{\times} \mc{C}_d \xrightarrow{\sim} (\mc{C}_d)_{\ell p/} .
			\]
		\end{enumerate}
	\end{lem}
	\begin{proof}
		For (i), it suffices to show that the left adjoint is defined on any object $C \in \mc{C}$. Let $C \to C'$ be a $q$-cocartesian arrow lying over the essentially unique arrow $q(C) \to d$. Then $C'$ corepresents $\ell(C)$.
		
		For (ii), the fact that $d$ is terminal implies that the left hand side is equivalent to $\mc{C}_{p/} \times_{\mc{C}} \mc{C}_d$. By the adjunction of (i), this is equivalent to $\mc{C}_{\ell p/}$.
	\end{proof}
	
	\begin{cor} \label{cor-underlying}
		Let $q : \mc{C}^\otimes \to \assoc$ be a monoidal $\infty$-category. Assume that $q$ is compatible with small colimits in the sense of \cite[3.1.1.19]{ha}. Let $\mc{I}$ be a small $\infty$-category, and let $F : \mc{I} \to \alg_{\assoc}(\mc{C})$ be a functor.
		\begin{enumerate}[label=(\roman*)]
			\item $\colim F$ exists in $\alg_{\assoc}(\mc{C})$.
			\item Applying Lemma~\ref{fill}(i) to the cocartesian fibration $\mc{C}^{\otimes}_{\on{act}} \to \assoc_{\on{act}}$ and the terminal object $\bb{1}\in \assoc_{\on{act}}$, we obtain a left adjoint $\ell : \mc{C}^{\otimes}_{\on{act}} \to \mc{C}^{\otimes}_{\on{act}, \bb{1}}$ to the embedding.
			Let $\wt{F} \in \alg_{\mc{I}^\amalg}(\mc{C})$ be the functor obtained from $F$ via Corollary~\ref{wr-cor}, and consider the composite functor defined as follows:
			\begin{cd}
				\mc{I}^\amalg_{\on{act}} \ar[r, "\wt{F}"] & \mc{C}^{\otimes}_{\on{act}} \ar[r, "\ell"] & \mc{C}^{\otimes}_{\on{act}, \bb{1}}
			\end{cd}
			Then $\colim (\ell \circ \wt{F}) \simeq \on{oblv}(\colim F)$.
			
			Here $\on{oblv} : \alg_{\assoc}(\mc{C}) \to \mc{C}$ is the forgetful functor.
		\end{enumerate}
	\end{cor}
	\begin{proof}
		By Corollary~\ref{wr-cor}, $\colim F$ is computed by applying the left adjoint of the restriction functor $\alg_{\assoc}(\oh) \xrightarrow{\on{Res}_{\pi'}} \alg_{\mc{I}^\amalg}(\oh)$ to $\wt{F} \in \alg_{\mc{I}^\amalg}(\mc{C})$. By~\cite[Cor.\ 3.1.3.5]{ha} and our hypothesis that $q$ is compatible with small colimits, this left adjoint exists and is given by the \emph{$q$-free $\assoc$-algebra generated by $\wt{F}$} in the sense of~\cite[Def.\ 3.1.3.1]{ha}. Applying condition ($*$) in that definition to the fibration $q : \mc{C}^\otimes \to \assoc$ and the maps $\mc{I}^\amalg \to \assoc \to \assoc$ of operads, with the object `$B$' taken to be $\bb{1} \in \assoc$ which is terminal in $\assoc_{\on{act}}$, and using that any operadic $q$-colimit diagram is a weak operadic $q$-colimit diagram (see~\cite[Def.\ 3.1.1.2]{ha}), we obtain that the underlying object $\on{oblv}(\colim F) \in \mc{C}^{\otimes}_{\bb{1}}$ corepresents
		\[
		\Big( (\mc{C}^\otimes_{\on{act}})_{\wt{F}/} \underset{(\assoc_{\on{act}})_{q\wt{F}/}}{\times}( \assoc_{\on{act}})_{p'/} \Big) \underset{\mc{C}^\otimes_{\on{act}}}{\times} \mc{C}^\otimes_{\on{act},\bb{1}}
		\]
		where the maps are defined as follows:
		\begin{cd}
			\mc{I}^{\amalg}_{\on{act}} \ar[r, "\wt{F}"] \ar[d, hookrightarrow] & \mc{C}^\otimes_{\on{act}} \ar[d, "q"] \\
			(\mc{I}^{\amalg}_{\on{act}})^{\triangleright} \ar[r, "p'"] & \assoc_{\on{act}}
		\end{cd}
		In this diagram, $p'$ sends the terminal vertex to $\bb{1} \in \assoc$. By Lemma~\ref{fill}(ii), this undercategory is equivalent to $\mc{C}_{\ell \wt{F} /}$, as desired.
	\end{proof}
	
	\subsubsection{Proof of Proposition~\ref{preamalg-prop}}
	Since $\emptyset \in \mc{P}_{I, \mr{fin}}$ is initial, we may think of $\mc{H}_\bullet$ as a functor to $\ms{Alg}(\cat)_{\mc{H}_\emptyset/}$. By~\cite[Cor.\ 3.4.1.7]{htt}, we have an equivalence of $\infty$-categories
	\[
		\ms{Alg}(\cat)_{\mc{H}_\emptyset/} \simeq \ms{Alg}({}_{\mc{H}_\emptyset}\ms{BMod}_{\mc{H}_\emptyset}(\cat))
	\]
	where the monoidal structure on the $\infty$-category of bimodules is the relative tensor product over $\mc{H}_\emptyset$. Apply Corollary~\ref{cor-underlying} with the following replacements: 
	\e{
		\mc{C} &:= {}_{\mc{H}_\emptyset}\ms{BMod}_{\mc{H}_\emptyset}(\cat) \\
		\mc{I} &:= \mc{P}_{I, \mr{fin}} \\ 
		F &:= \mc{H}_\bullet
	} 
	The fiber $\infty$-category $\mc{C}^\otimes_{\mr{act}, \bb{1}}$ is equivalent to $\mc{C}$, and the functor 
	\[
		\ell \circ \wt{F} : \mc{I}^\amalg_{\mr{act}} \to \mc{C}
	\]
	gives the desired functor $\wh{\mc{H}}_\bullet$ upon post-composing by the forgetful functor 
	\[
		{}_{\mc{H}_\emptyset}\ms{BMod}_{\mc{H}_\emptyset}(\cat) \to \cat
	\]
	In the statement of Proposition~\ref{preamalg-prop}, the first two desired properties follow by unwinding the operadic colimits in the proof of Corollary~\ref{cor-underlying}, and the third desired property is Corollary~\ref{cor-underlying}(ii). \hfill $\qed$ 
	
	\needspace{1.5in} \subsection{Generalized Demazure varieties} \label{ss-interpret}
	
	\subsubsectiona \label{interpret-def} 
	In~\ref{firstdem}, we defined a functor $D^{\ms{d}} : \ms{Rig}^{\ms{d}} \to \ms{Varieties}$ whose values are twisted products of Schubert varieties. We review and generalize this definition, using different notations which are better suited to the current setting. 
	
	For every element $w \in W$, let $\mb{P}_w \subset \LG$ be the closure of the $w$ Bruhat cell, so that $\Fl_w := \mb{P}_w / \mb{I} \subset \Fl$ is the $w$ Schubert variety. Also, let $\mc{H}_w \subset \mc{H}_I$ be the full subcategory consisting of sheaves supported on $\Fl_w$. 
	
	Let $T$ be a maximal torus of $G$, and let $\mb{U} \subset \mb{I}$ be the pro-unipotent radical, so we have an exact sequence $\mb{U} \to \mb{I} \to T$. For any $\bm{w} = ([1,n], w, J) \in \ms{Rig}^{\ms{d}}$, we define the following: 
	\e{
		\wt{\Fl}_{\bm{w}} &:= \mathbf{P}_{w_1} \overset{\mathbf{I}}{\times} \mathbf{P}_{w_2} \overset{\mathbf{I}}{\times} \cdots \overset{\mathbf{I}}{\times} \mathbf{P}_{w_n} / \mathbf{I} \\
		\ol{\Fl}_{\bm{w}} &:= \big(\mathbf{P}_{w_1}/\mathbf{U}\big) \overset{T}{\times} \big(\mathbf{U} \backslash \mathbf{P}_{w_2} / \mathbf{U}\big) \overset{T}{\times} \cdots \overset{T}{\times} \big(\mathbf{U} \backslash \mathbf{P}_{w_n} / \mathbf{I}\big)
	}
	Note that $\wt{\Fl}_{\bm{w}}$ is a twisted product of Schubert varieties and was previously denoted $D^{\ms{d}}([1,n], w, J)$. We call it a \emph{generalized Demazure variety}. The multiplication map 
	\[
		m : \wt{\Fl}_{\bm{w}} \to \Fl
	\]
	will play an important role. In contrast, $\ol{\Fl}_{\bm{w}}$ is a stack with pro-unipotent inertia groups, and it does not admit such a map. 
	
	The Bruhat stratification of $\mb{P}_{w_i}$ induces stratifications of $\wt{\Fl}_{\bm{w}}$ and $\ol{\Fl}_{\bm{w}}$, whose cells are indexed by tuples $(w'_1, \ldots, w'_n)$ such that $w'_i \preceq w_i$. Note that these are in bijection with $\ms{Emb}^{\ms{d}}$-maps to $\bm{w}$. Let
	\[
		\Shv'(\mb{I} \backslash \wt{\Fl}_{\bm{w}}) \subset \Shv(\mb{I} \backslash \wt{\Fl}_{\bm{w}}) 
	\]
	be the full subcategory consisting of $\Shv$-modules whose $*$-restriction to each cell is constant. We have an adjunction 
	\[
		\Shv'(\mb{I} \backslash \wt{\Fl}_{\bm{w}}) \underset{a_*}{\overset{a^*}{\rightleftarrows}} \Shv(\mb{I} \backslash \wt{\Fl}_{\bm{w}}) 
	\]
	where the left adjoint $a^*$ is the embedding and the right adjoint $a_*$ can be informally described as `twisted $*$-averaging with respect to $\mb{U}$-actions.' 
	
	\begin{lem} \label{interpret-obj}
		For every $\bm{w} = ([1, n], w, J)$ as above, we have canonical isomorphisms
		\[
			\mc{H}_{w_1} \underset{\mc{H}_\emptyset}{\otimes} \cdots \underset{\mc{H}_\emptyset}{\otimes} \mc{H}_{w_n} \simeq \Shv(\mathbf{I} \backslash \ol{\Fl}_{\bm{w}}) \simeq \Shv'(\mb{I} \backslash \wt{\Fl}_{\bm{w}})
		\]
	\end{lem}
	\begin{proof}
		For notational clarity, let us assume $n = 2$. By definition, we have 
		\[
			\mc{H}_{w_1} \underset{\mc{H}_{\emptyset}}{\otimes} \mc{H}_{w_2} := \Shv(\mb{I} \backslash \mb{P}_{w_1} / \mb{I}) \underset{\Shv(\mb{I} \backslash \mb{I} / \mb{I})}{\otimes} \Shv(\mb{I} \backslash \mb{P}_{w_2} / \mb{I})
		\]
		We will simplify this tensor product. The convolution monoidal structure on $\Shv(\mb{I} \backslash \mb{I} / \mb{I})$ canonically identifies with the ordinary monoidal structure given by tensor product of $\D$-modules. Furthermore, we have 
		\[
		\Shv(\mb{I} \backslash \mb{I} / \mb{I}) \simeq \Shv(\mr{pt}/\mb{I}) \simeq \Shv(\mr{pt}/T), 
		\]
		where the second equivalence follows from the fact that $\mb{U}$ is pro-unipotent, hence contractible. Under these identifications, the right $\mc{H}_{\emptyset}$-module structure of $\mc{H}_{w_1}$ agrees with the right $\Shv(\mr{pt}/T)$-module structure obtained by pulling back $\D$-modules along the projection map 
		\[
		p_1 : \mb{I} \backslash \mb{P}_{w_1} / \mb{I} \to \mr{pt} / T
		\]
		which comes from the right $T$-action on $\mb{I} \backslash \mb{P}_{w_1} / \mb{U}$. Similarly, the left $\mc{H}_{\emptyset}$-module structure of $\mc{H}_{w_2}$ agrees with the left $\Shv(\mr{pt}/T)$-module structure obtained by pulling back $\D$-modules along the projection map 
		\[
		p_2 : \mb{I} \backslash \mb{P}_{w_2} / \mb{I} \to T \backslash \mr{pt}
		\]
		which comes from the left $T$-action on $\mb{U} \backslash \mb{P}_{w_2} / \mb{I}$. In summary, we have 	
		\[
			\Shv(\mb{I} \backslash \mb{P}_{w_1} / \mb{I}) \underset{\Shv(\mb{I} \backslash \mb{I} / \mb{I})}{\otimes} \Shv(\mb{I} \backslash \mb{P}_{w_2} / \mb{I}) \simeq \Shv(\mb{I} \backslash \mb{P}_{w_1} / \mb{I}) \underset{\Shv(\pt / T)}{\otimes} \Shv(\mb{I} \backslash \mb{P}_{w_2} / \mb{I})
		\]
		
		Next, we will construct the first isomorphism in the statement of the lemma. For a suitably nice pullback square 
		\begin{cd}
			\mc{X} \times_{\mc{Z}} \mc{Y} \ar[r] \ar[d] & \mc{Y} \ar[d] \\
			\mc{X} \ar[r] & \mc{Z}
		\end{cd}
		of stacks,~\cite[Thm.\ 1.15]{ben} shows that the exact functor $\Shv(\mc{X}) \otimes_{\Shv(\mc{Z})} \Shv(\mc{Y}) \to \Shv(\mc{X} \times_{\mc{Z}} \mc{Y})$ is fully faithful. In order to apply the theorem, we first consider this pullback square of stacks: 
		\begin{cd}
			\mb{I} \backslash \ol{\Fl}_{\bm{w}} \ar[r] \ar[d] & \mb{I} \backslash \mb{P}_{w_2} / \mb{I} \ar[d, "p_2"] \\
			\mb{I} \backslash \mb{P}_{w_1} / \mb{I} \ar[r, "p_1"] & \mr{pt} / T
		\end{cd}
		Let $\mb{U}' \subset \mb{U}$ be a subgroup such that $\mb{U}'$ is normal in $\mb{I}$, $\mb{U}'$ acts trivially on $\mb{I} \backslash \mb{P}_{w_1}$ from the right and on $\mb{P}_{w_2} / \mb{I}$ from the left, and the quotient $\mb{U}'' := \mb{U} / \mb{U}'$ is finite type. Since $\mb{U}$ is pro-unipotent, the $\D$-module categories in question are unchanged if we replace the above diagram by the following one: 
		\begin{cd}
			\mb{I} \backslash \mb{P}_{w_1} / \mb{U}'' \overset{T}{\times} \mb{U}'' \backslash \mb{P}_{w_2} / \mb{I} \ar[r] \ar[d] & (\mb{I} / \mb{U}') \backslash \mb{P}_{w_2} / \mb{I} \ar[d, "p_2"] \\
			\mb{I} \backslash \mb{P}_{w_1} / (\mb{I} / \mb{U}') \ar[r, "p_1"] & \mr{pt} / T
		\end{cd}
		Finally,~\cite[Thm.\ 1.15]{ben} applies to this diagram because the morphisms $p_1$ and $p_2$ are \emph{safe}. This assertion follows from the definition of `safe morphism' (see \cite[Def.\ 10.2.2]{finiteness}) since $p_1$ and $p_2$ are quasicompact morphisms of algebraic stacks, such that the geometric fibers $\mb{I} \backslash \mb{P}_{w_1} / \mb{U}''$ and $\mb{U}'' \backslash \mb{P}_{w_2} / \mb{I}$ satisfy the property that every inertia group is unipotent. Hence, the theorem implies that 
		\[
			\Shv(\mb{I} \backslash \mb{P}_{w_1} / \mb{I}) \underset{\Shv(\pt / T)}{\otimes} \Shv(\mb{I} \backslash \mb{P}_{w_2} / \mb{I}) \to \Shv(\mb{I} \backslash \ol{\Fl}_{\bm{w}})
		\]
		is fully faithful.  
		
		To show essential surjectivity, note that $\Shv(\mb{I} \backslash \ol{\Fl}_{\bm{w}})$ is generated by the $!$-pushforwards of constant sheaves on the open cells. These objects clearly lie in the image of the aforementioned functor. Since the functor is exact and fully faithful, this shows that it is essentially surjective. We have now constructed the first equivalence in the lemma. 
		
		For the second equivalence, which states that $\Shv(\mathbf{I} \backslash \ol{\Fl}_{\bm{w}}) \simeq \Shv'(\mb{I} \backslash \wt{\Fl}_{\bm{w}})$, note that the $*$-pullback functor 
		\[
			\Shv(\mathbf{I} \backslash \ol{\Fl}_{\bm{w}}) \to \Shv(\mb{I} \backslash \wt{\Fl}_{\bm{w}})
		\]
		is exact and fully faithful (since $\mb{U}$ is contractible). This functor is an equivalence onto the full subcategory $\Shv'(\mb{I} \backslash \wt{\Fl}_{\bm{w}})$ because this $\infty$-category and $\Shv(\mathbf{I} \backslash \ol{\Fl}_{\bm{w}})$ are both generated by $!$-pushforwards of constant sheaves on open cells. 
	\end{proof}
	
	\begin{cor} \label{interpret-cor} 
		There exists a functor 
		\[
		\wt{\mc{H}}_\bullet : \ms{Rig}^{\ms{d}} \to \cat
		\]
		with the following properties: 
		\begin{itemize}
			\item For every $\bm{w} = ([1, n], w, J) \in \ms{Rig}^{\ms{d}}$, we have 
			\[
				\wt{\mc{H}}_{\bm{w}} = \mc{H}_{w_1} \underset{\mc{H}_\emptyset}{\otimes} \cdots \underset{\mc{H}_\emptyset}{\otimes} \mc{H}_{w_n}
			\]
			\item $\wt{\mc{H}}_\bullet$ sends a $\ms{Rig}^{\ms{d}}$-morphism $([1, n], w, J) \xra{\psi : [1,n] \to [1,n']} ([1, n'], w', J')$ to the functor described as follows. For each $j \in [1, n']$, write $\psi^{-1}(j) = [a_j, b_j]$, and consider the functor  
			\[
			\mc{H}_{w_{a_j}} \underset{\mc{H}_\emptyset}{\otimes} \mc{H}_{w_{a_j+1}} \underset{\mc{H}_\emptyset}{\otimes} \cdots \underset{\mc{H}_\emptyset}{\otimes} \mc{H}_{w_{b_j}} \to \mc{H}_{w'_j}
			\]
			induced by the convolution monoidal structure of $\mc{H}_{I}$. Then tensor over $j \in [1, n']$. 
			\item There is an equivalence of DG-categories 
			\[
				\colim_{\ms{Rig}^{\ms{d}}} \wt{\mc{H}}_\bullet \xra{\sim} \mr{oblv}(\colim_{\mc{P}_{I, \mr{fin}}} \mc{H}_\bullet). 
			\]
		\end{itemize}
	\end{cor}
	In the last bullet, `oblv' refers to forgetting the monoidal structure. 
	\begin{proof}
		First, there is an adjunction
		\[
		\ms{Rig}^{\ms{d}} \overset{L}{\underset{R}{\rightleftarrows}} (\mc{P}_{I, \mr{fin}})^{\amalg}_{\on{act}}
		\]
		where the left adjoint $L$ forgets the sequence $w$, and the right adjoint sends $([1, n], J) \mapsto ([1, n], w^{\ms{max}}, J)$, where $w^{\ms{max}}_i$ is the longest element of $W_{J_i}$. Since $R$ is fully faithful, $L \circ R \simeq \id$. 
		
		Define $\wt{\mc{H}}_\bullet$ to be the subfunctor of $\wh{\mc{H}}_\bullet \circ L$ whose behavior on objects is as follows: for each $\bm{w} = ([1, n], w, J) \in \ms{Rig}^{\ms{d}}$, we have 
		\[
			\wt{\mc{H}}_{\bm{w}} := \mc{H}_{w_1} \underset{\mc{H}_\emptyset}{\otimes} \cdots \underset{\mc{H}_\emptyset}{\otimes} \mc{H}_{w_n} \hra \mc{H}_{J_1} \underset{\mc{H}_\emptyset}{\otimes} \cdots \underset{\mc{H}_\emptyset}{\otimes} \mc{H}_{J_n} = \wh{\mc{H}}_{L(\bm{w})}. 
		\]
		This map is fully faithful because, under Lemma~\ref{interpret-obj}, it corresponds to $*$-pushforward of $\D$-modules along a closed embedding. The behavior of $\wt{\mc{H}}_\bullet$ on arrows and higher simplices is completely determined by $\wh{\mc{H}}_\bullet \circ L$. To check that $\wt{\mc{H}}_\bullet$ is well-defined, we only need to check that, for each map $\bm{w} \xra{\psi : [1,n'] \to [1,n]} \bm{w}'$, the map 
		\[
			\wh{\mc{H}}_{L(\psi)} : \wh{\mc{H}}_{L(\bm{w})} \to \wh{\mc{H}}_{L(\bm{w}')}
		\]
		sends the full subcategory $\wt{\mc{H}}_{\bm{w}}$ into $\wt{\mc{H}}_{\bm{w}'}$. This follows from the condition
		\[
			\big(\text{Demazure product of $(\psi^{-1}(j), w)$}\big) \preceq w'_j
		\]
		which applies to morphisms in $\ms{Rig}^{\ms{d}}$, together with the fact that the Demazure product governs the image of Schubert varieties under convolution. 
		
		The first two bullets are true by construction. For the third bullet, note that 
		\e{
			\colim_{\ms{Rig}^{\ms{d}}} \wt{\mc{H}}_\bullet &\simeq \colim_{(\mc{P}_{I, \mr{fin}})^{\amalg}_{\on{act}}} \wt{\mc{H}}_\bullet \circ R \\
			&\simeq \colim_{(\mc{P}_{I, \mr{fin}})^{\amalg}_{\on{act}}} \wh{\mc{H}}_\bullet \\
			&\simeq \mr{oblv}(\colim_{\mc{P}_{I, \mr{fin}}} \mc{H}_\bullet), 
		} 
		where the first equivalence follows from the fact that $R$ is final (because it is a right adjoint), the second follows from $\wt{\mc{H}}_\bullet \circ R \simeq \wh{\mc{H}}_\bullet$, and the third is part of Proposition~\ref{preamalg-prop}. 
	\end{proof}	
	
	\subsubsectiona \label{interpret-arr} 
	The following lemma can be used to describe the behavior of $\wt{\mc{H}}_\bullet$ on arrows. 
	\begin{lem*} 
		For any $\bm{w} = ([1, n], w, J) \in \ms{Rig}^{\ms{d}}$ and any $w' \in W_I$ satisfying 
		\[
			\big(\text{Demazure product of $([1,n], w)$}\big) \preceq w'
		\]
		we have a commutative diagram 
		\begin{cd}
			\mc{H}_{w_1} \underset{\mc{H}_\emptyset}{\otimes} \cdots \underset{\mc{H}_\emptyset}{\otimes} \mc{H}_{w_n} \ar[rr] \ar[dash]{d}[rotate=90, anchor=north]{\sim} & & \mc{H}_{w'} \ar[dash]{d}[rotate=90, anchor=north]{\sim} \\
			\Shv'(\mb{I} \backslash \wt{\Fl}_{\bm{w}}) \ar[r, "a^*"] & \Shv(\mb{I} \backslash \wt{\Fl}_{\bm{w}}) \ar[r, "m_!"] & \Shv(\mb{I} \backslash  \Fl_{w'})
		\end{cd}
		The upper horizontal map is induced from the monoidal structure of $\mc{H}_I$, the vertical maps come from Lemma~\ref{interpret-obj}, and the lower horizontal maps are from~\ref{interpret-def}. 
	\end{lem*}
	\begin{proof}
		For notational clarity, assume $n = 2$. By the construction of the convolution tensor product on $\mc{H}_I$, we have a commutative diagram 
		\begin{cd}
			\mc{H}_{w_1} \otimes \mc{H}_{w_2} \ar[rr] \ar[dash]{d}[rotate=90, anchor=north]{\sim} & & \mc{H}_{w'} \ar[dash]{d}[rotate=90, anchor=north]{\sim} \\
			\Shv\big(\big(\mb{I} \backslash \bm{P}_{w_1} / \mb{I}\big) \times \big(\mb{I} \backslash \bm{P}_{w_2} / \mb{I}\big)\big) \ar[r, "*\text{-pull}"] & \Shv(\mb{I} \backslash \wt{\Fl}_{\bm{w}}) \ar[r, "m_!"] & \Shv(\mb{I} \backslash \Fl_{w'})
		\end{cd}
		where the upper horizontal map is induced by the tensor product. By factoring the projection map $\mb{I} \backslash \wt{\Fl}_{\bm{w}} \to \big(\mb{I} \backslash \bm{P}_{w_1} / \mb{I}\big) \times \big(\mb{I} \backslash \bm{P}_{w_2} / \mb{I}\big)$ as a composite 
		\begin{cd}
			\mb{I} \backslash \wt{\Fl}_{\bm{w}} \ar[r] & \mb{I} \backslash \ol{\Fl}_{\bm{w}} \ar[r] & 
			\big(\mb{I} \backslash \bm{P}_{w_1} / \mb{I}\big) \times \big(\mb{I} \backslash \bm{P}_{w_1} / \mb{I}\big)
		\end{cd}
		we can enlarge the previous diagram into a solid diagram as follows: 
		\begin{cd}
			\mc{H}_{w_1} \otimes \mc{H}_{w_2} \ar[r] \ar[dash]{d}[rotate=90, anchor=north]{\sim} & \mc{H}_{w_1} \underset{\mc{H}_\emptyset}{\otimes}  \mc{H}_{w_2} \ar[rr] \ar[d, dashed] & & \mc{H}_{w'} \ar[dash]{d}[rotate=90, anchor=north]{\sim} \\
			\Shv\big(\big(\mb{I} \backslash \bm{P}_{w_1} / \mb{I}\big) \times \big(\mb{I} \backslash \bm{P}_{w_1} / \mb{I}\big)\big) \ar[r, "*\text{-pull}"] & \Shv(\mb{I} \backslash \ol{\Fl}_{\bm{w}}) \ar[r, "*\text{-pull}"] & \Shv(\mb{I} \backslash \wt{\Fl}_{\bm{w}}) \ar[r, "m_!"] & \Shv(\mb{I} \backslash \Fl_{w'})
		\end{cd}
		The universal property of relative tensor product yields a dashed arrow which makes the diagram commute, and tracing through definitions shows that the dashed arrow identifies with the first isomorphism in Lemma~\ref{interpret-obj}. Now the desired statement follows from the second isomorphism in Lemma~\ref{interpret-obj}. 
	\end{proof}
	
	\mysubsec{Main theorem} \label{ss-checkbis}
	
	\subsubsectiona \label{blowup} 
	For any $\bm{w} \in \ms{Rig}^{\ms{d}}$, let $\wt{\Fl}_{\bm{w}}^\circ, \partial \wt{\Fl}_{\bm{w}} \hra \wt{\Fl}_{\bm{w}}$ be the open cell and its closed complement. 
	\begin{lem*} 
		Suppose that $\bm{w}_1 \to \bm{w}_2$ is a basic level map in $\ms{Rig}^{\ms{d}}$. Then $m : \wt{\Fl}_{\bm{w}_1} \to \wt{\Fl}_{\bm{w}_2}$ is an isomorphism on the open strata, and the commutative diagram 
		\begin{cd}
			\Shv'(\mb{I} \backslash \partial \wt{\Fl}_{\bm{w}_1}) \ar[r, hookrightarrow] \ar[d] & \Shv'(\mb{I} \backslash \wt{\Fl}_{\bm{w}_1}) \ar[d, "m_!"] \\
			\Shv'(\mb{I} \backslash \partial \wt{\Fl}_{\bm{w}_2}) \ar[r, hookrightarrow] & \Shv'(\mb{I} \backslash \wt{\Fl}_{\bm{w}_2})
		\end{cd}
		is cocartesian. 
	\end{lem*}
	\begin{proof}
		Recall that basic level maps in $\ms{Rig}^{\ms{d}}$ correspond to level maps in $\ms{Dom}^{\ms{c}}$ (\ref{rigidcat-deform}), and a map in $\ms{Dom}^{\ms{c}}$ is level if and only if its Coxeter products involve reduced sequences (\ref{good}). This implies that $m$ is an isomorphism on the open strata. 
		
		Next, we prove that the above diagram is cocartesian. We will use the diagram 
		\begin{cd}
			\partial \wt{\Fl}_{\bm{w}_1} \ar[r, hookrightarrow, "i"] \ar[d, "\dot{m}"] & \wt{\Fl}_{\bm{w}_1} \ar[r, hookleftarrow, "j"] \ar[d, "m"] & \wt{\Fl}^\circ_{\bm{w}_1} \ar{d}{\ddot{m}}[swap, rotate = 90, anchor = south]{\sim} \\
			\partial \wt{\Fl}_{\bm{w}_2} \ar[r, hookrightarrow, "\imath"] & \wt{\Fl}_{\bm{w}_2} \ar[r, hookleftarrow, "\jmath"] &  \wt{\Fl}^\circ_{\bm{w}_2}
		\end{cd}
		
		Since the Schubert stratification is Whitney, the functors $j_!$ and $j_*$ preserve the subcategories of $\D$-modules denoted by $\Shv'(-)$, so we have adjoint triples 
		\begin{cd}[column sep = 0.7in] 
			\Shv'(\mb{I} \backslash \wt{\Fl}_{\bm{w}_1}^\circ) \ar[r, bend left = 20, "j_!" description] \ar[r, bend right = 20, "j_*" description] & \Shv'(\mb{I} \backslash \wt{\Fl}_{\bm{w}_1}) \ar[l, "j^*" description] \ar[r, bend left = 20, "i^*" description] \ar[r, bend right = 20, "i^!" description] & \Shv'(\mb{I} \backslash \partial\wt{\Fl}_{\bm{w}_1}) \ar[l, "i_*" description] 
		\end{cd}
		By the main theorem of~\cite{strat} (see also Equation (5) therein), this yields a `recollement' description of $\Shv'(\mb{I} \backslash \wt{\Fl}_{\bm{w}_1})$, namely
		\[
			\Shv'(\mb{I} \backslash \wt{\Fl}_{\bm{w}_1}) \simeq {\lim}^{\ms{r.lax}} \Big( \begin{tikzcd} 
				\Shv'(\mb{I} \backslash \wt{\Fl}_{\bm{w}_1}^\circ) \ar[r, "i^*j_*"] & \Shv'(\mb{I} \backslash \partial \wt{\Fl}_{\bm{w}_1}) \end{tikzcd} \Big). 
		\]
		The notation $\lim^{\ms{r.lax}}$ means the right-lax limit in the $(\infty,2)$-category of $(\infty,1)$-categories. The Verdier dual statement also holds, and it is the one we will actually use: 
		\[
		\Shv'(\mb{I} \backslash \wt{\Fl}_{\bm{w}_1}) \simeq {\lim}^{\ms{l.lax}} \Big( \begin{tikzcd} 
		\Shv'(\mb{I} \backslash \wt{\Fl}_{\bm{w}_1}^\circ) \ar[r, "i^!j_!"] & \Shv'(\mb{I} \backslash \partial \wt{\Fl}_{\bm{w}_1}) \end{tikzcd} \Big). 
		\]
		As remarked in~\cite[0.2.1]{strat}, the right or left lax limit can be concretely described as a comma $\infty$-category. That is, the right hand side of the preceding equation can be described as the $\infty$-category of triples 
		\[
			\begin{pmatrix} \begin{array}{rcl}
			\mc{F}  &\in&\Shv'(\mb{I} \backslash \wt{\Fl}_{\bm{w}_1}^\circ)\\ 
			\mc{G} &\in& \Shv'(\mb{I} \backslash \partial \wt{\Fl}_{\bm{w}_1}) \\
			\varphi &:& i^!j_!\mc{F} \to \mc{G} 
			\end{array} \end{pmatrix}. 
		\]
		Analogous statements hold when $\bm{w}_1$ is replaced by $\bm{w}_2$. 
		
		By~\cite[Cor.\ 5.5.3.4]{htt}, it is equivalent to show that the diagram obtained by passing to right adjoints is cartesian. Thanks to the previous paragraph, we may rewrite the right adjoint diagram as follows: 
		\begin{cd}
			{\lim}^{\ms{l.lax}} \Big( 
				\Shv'(\mb{I} \backslash \wt{\Fl}_{\bm{w}_2}^\circ) \xra{\imath^!\jmath_!} \Shv'(\mb{I} \backslash \partial \wt{\Fl}_{\bm{w}_2}) \Big) \ar[r] \ar[d] & \Shv'(\mb{I} \backslash \partial \wt{\Fl}_{\bm{w}_2}) \ar[d, "a_*\dot{m}^!"] \\
			{\lim}^{\ms{l.lax}} \Big( 
				\Shv'(\mb{I} \backslash \wt{\Fl}_{\bm{w}_1}^\circ) \xra{i^!j_!}\Shv'(\mb{I} \backslash \partial \wt{\Fl}_{\bm{w}_1})  \Big) \ar[r] & \Shv'(\mb{I} \backslash \partial \wt{\Fl}_{\bm{w}_1}) 
		\end{cd}
		(The `averaging' functor $a_*$ was defined in~\ref{interpret-def}.) Abstractly, the left vertical map is induced by taking left-lax limits along the rows of this lax-commutative diagram: 
		\begin{cd}
			\Shv'(\mb{I}\backslash\wt{\Fl}^\circ_{\bm{w}_2})\ar[r, "\imath^!\jmath_!"] \ar{d}{a_*\ddot{m}^!}[swap, anchor=south, rotate=90]{\sim} & \Shv'(\mb{I}\backslash\partial \wt{\Fl}_{\bm{w}_2}) \ar[d, "a_*\dot{m}^!"] \\
			\Shv'(\mb{I}\backslash\wt{\Fl}^\circ_{\bm{w}_1})\ar[r, "i^!j_!"] \ar[draw=none]{ru}[anchor=center, rotate=135]{\Downarrow} & \Shv'(\mb{I}\backslash\partial \wt{\Fl}_{\bm{w}_1})
		\end{cd}
		The natural transformation (not isomorphism!)\ witnessing this lax-commutativity is 
		\[
			i^!j_!a_*\ddot{m}^! \simeq a_* i^! j_! \ddot{m}^! 
			\xra{\eta} a_* i^! m^! \jmath_! 
			\simeq a_* \dot{m}^! \imath^! \jmath_!
		\]
		This uses the base-change natural transformation $j_! \ddot{m}^! \xra{\eta} m^! \jmath_!$. Concretely, the left vertical map sends 
		\[
			\begin{pmatrix}
			\begin{array}{rcl}
			\mc{F} &\in& \Shv'(\mb{I}\backslash \wt{\Fl}_{\bm{w}_2}^\circ)\\
			\mc{G} &\in& \Shv'(\mb{I}\backslash \partial\wt{\Fl}_{\bm{w}_2}) 
			\end{array} \\
			\imath^!\jmath_! \mc{F} \xra{\varphi} \mc{G}
			\end{pmatrix} \mapsto \begin{pmatrix}
			\begin{array}{rcl}
			a_*\ddot{m}^!\mc{F} &\in& \Shv'(\mb{I}\backslash \wt{\Fl}_{\bm{w}_1}^\circ) \\
			a_*\dot{m}^!\mc{G} &\in& \Shv'(\mb{I}\backslash \partial\wt{\Fl}_{\bm{w}_1})
			\end{array} \\
			i^!j_! a_*\ddot{m}^! \mc{F} \to a_* \dot{m}^! \imath^! \jmath_! \mc{F} \xra{a_*\dot{m}^!(\varphi)} a_*\dot{m}^! \mc{G}
			\end{pmatrix},  
		\]
		where the unlabeled arrow comes from the lax-commutativity natural transformation. 
		
		In the right adjoint diagram, the fibered product of the lower-right cospan can be concretely described as follows: it consists of triples 
		\[
			\begin{pmatrix}
			\begin{array}{rcl}
			\mc{F} &\in& \Shv'(\mb{I} \backslash \wt{\Fl}_{\bm{w}_1}^\circ) \\
			\mc{G} &\in& \Shv'(\mb{I} \backslash \partial \wt{\Fl}_{\bm{w}_2}) 
			\end{array} \\
			i^!j_! \mc{F} \xra{\varphi} a_*\dot{m}^! \mc{G}
		\end{pmatrix},  \text{ or equivalently } 	\begin{pmatrix}
		\begin{array}{rcl}
		\mc{F} &\in& \Shv'(\mb{I} \backslash \wt{\Fl}_{\bm{w}_1}^\circ) \\
		\mc{G} &\in& \Shv'(\mb{I} \backslash \partial \wt{\Fl}_{\bm{w}_2}) 
		\end{array} \\
		\dot{m}_* i^!j_! \mc{F} \xra{\varphi} \mc{G}
		\end{pmatrix}, 
		\]
		where the equivalence follows from the adjunction $(\dot{m}_*, \dot{m}^!)$ and the observation that $j_!$ and $\dot{m}_*$ preserve the subcategory denoted by $\Shv'(-)$. The map from the upper-left vertex to the fibered product can be described as follows: 
		\[
			\begin{pmatrix}
			\begin{array}{rcl}
			\mc{F} &\in& \Shv'(\mb{I}\backslash \wt{\Fl}_{\bm{w}_2}^\circ)\\
			\mc{G} &\in& \Shv'(\mb{I}\backslash \partial\wt{\Fl}_{\bm{w}_2}) 
			\end{array} \\
			\imath^!\jmath_! \mc{F} \xra{\varphi} \mc{G}
			\end{pmatrix} \mapsto \begin{pmatrix}
			\begin{array}{rcl}
			a_*\ddot{m}^!\mc{F} &\in& \Shv'(\mb{I}\backslash \wt{\Fl}_{\bm{w}_1}^\circ) \\
			\mc{G} &\in& \Shv'(\mb{I}\backslash \partial\wt{\Fl}_{\bm{w}_2}) 
			\end{array} \\
			i^!j_! a_*\ddot{m}^! \mc{F} \to a_* \dot{m}^! \imath^! j_! \mc{F} \xra{a_*\dot{m}^!(\varphi)} a_*\dot{m}^! \mc{G}
			\end{pmatrix}. 
		\]
		To prove that the diagram is cartesian, we write down the inverse functor: 
		\[
			\begin{pmatrix}
			\begin{array}{rcl}
			\mc{F} &\in& \Shv'(\mb{I} \backslash \wt{\Fl}_{\bm{w}_1}^\circ) \\
			\mc{G} &\in& \Shv'(\mb{I} \backslash \partial \wt{\Fl}_{\bm{w}_2}) 
			\end{array} \\
			\dot{m}_*i^!j_! \mc{F} \xra{\varphi} \mc{G}
			\end{pmatrix} \mapsto \begin{pmatrix}
			\begin{array}{rcl}
				\ddot{m}_*\mc{F} &\in& \Shv'(\mb{I} \backslash \wt{\Fl}_{\bm{w}_2}^\circ) \\
				\mc{G} &\in& \Shv'(\mb{I} \backslash \partial \wt{\Fl}_{\bm{w}_2}) 
			\end{array} \\
				\imath^!\jmath_! \ddot{m}_*\mc{F} \simeq \imath^! m_* j_! \mc{F} \simeq \dot{m}_* i^!j_! \mc{F} \xra{\varphi} \mc{G}
			\end{pmatrix}
		\]
		The definition of the map on the right hand side uses $\ddot{m}_* \simeq \ddot{m}_!$ and $\imath^!m_* \simeq \dot{m}_*i^!$. We leave to the reader the straightforward but tedious exercise of checking that this does indeed give an inverse functor. 		
	\end{proof}
	
	\begin{cor} \label{blowup-cor} 
		The functor $\wt{\mc{H}}_\bullet : \ms{Rig}^{\ms{d}} \to \cat$ satisfies~\ref{smooth}(3). 
	\end{cor}
	\begin{proof}
		First, fix $\bm{w} \in \ms{Rig}^{\ms{d}}$ and write $j := r(\bm{w})$. We will prove that there is a dashed equivalence making the following diagram commute. 
		\begin{cd}
			\displaystyle\colim_{\la \ms{Emb}^{\ms{d}}_{<j} \xra{\ms{Emb}^{\ms{d}}} \bm{w} \ra} \wt{\mc{H}}_\bullet \ar[r, dashed, "\sim"] \ar[d] & \Shv'(\mb{I} \backslash \partial \wt{\Fl}_{\bm{w}}) \ar[d, hookrightarrow] \\
			\wt{\mc{H}}_{\bm{w}} \ar{r}{\sim}[swap]{\text{Lemma~\ref{interpret-obj}}} & \Shv'(\mb{I} \backslash \wt{\Fl}_{\bm{w}})
		\end{cd}
		To construct the dashed equivalence, first use Lemma~\ref{interpret-obj} to write 
		\[
			\colim_{\la \ms{Emb}^{\ms{d}}_{<j} \xra{\ms{Emb}^{\ms{d}}} \bm{w} \ra} \wt{\mc{H}}_\bullet \simeq \colim_{(\bm{w}' \hra \bm{w}) \in \la \ms{Emb}^{\ms{d}}_{<j} \xra{\ms{Emb}^{\ms{d}}} \bm{w} \ra} \Shv'(\bm{I} \backslash \wt{\Fl}_{\bm{w}'}), 
		\]
		and note that the colimit ranges over all non-identity $\ms{Emb}^{\ms{d}}$-maps to $\bm{w}$. Then the dashed equivalence will be the obvious map  
		\[ 
			\colim_{\bm{w}' \hra \bm{w}} \Shv'(\bm{I} \backslash \wt{\Fl}_{\bm{w}'}) \to \Shv'(\bm{I} \backslash \partial\wt{\Fl}_{\bm{w}}). 
		\]
		To show that it is an equivalence, it is convenient to pass to its right adjoint map 
		\[
			\Shv'(\bm{I} \backslash \partial\wt{\Fl}_{\bm{w}}) \to \lim_{\bm{w}' \hra \bm{w}} \Shv'(\bm{I} \backslash \wt{\Fl}_{\bm{w}'}), 
		\]
		where the limit diagram is defined via $!$-pullbacks. (This characterization of the right adjoint map follows from~\cite[Cor.\ 5.5.3.4]{htt}.) The method of Lemma~\ref{blowup} implies that this right adjoint map is an equivalence. Indeed, using the main theorem of~\cite{strat}, one can describe $\Shv'(\bm{I} \backslash \partial\wt{\Fl}_{\bm{w}})$ as the lax limit of the $\infty$-categories of sheaves on each cell, and similarly for the $\infty$-categories appearing in the limit. 
		
		Now let $\bm{w} \to \bm{w}'$ be a basic level morphism in $\ms{Rig}^{\ms{d}}$, which gives a diagram as in~\ref{smooth}(3). The previous paragraph implies that the horizontal arrows in~\ref{smooth}(3) are fully faithful, so the commutativity datum in~\ref{smooth}(3) is determined by the right vertical map  $\wt{\mc{H}}_{\bm{w}} \to \wt{\mc{H}}_{\bm{w}'}$. We may then use the equivalence in the the previous paragraph to identify the diagram of~\ref{smooth}(3) with that of Lemma~\ref{blowup}, which we have already shown is cocartesian. 
	\end{proof}
	
	\begin{rmk*}
		Here is a slightly different proof that the right adjoint map (from above) is an equivalence. First note that 
		\[
		\Shv(\bm{I} \backslash \partial\wt{\Fl}_{\bm{w}}) \to \lim_{\bm{w}' \hra \bm{w}} \Shv(\bm{I} \backslash \wt{\Fl}_{\bm{w}'})
		\]
		is an equivalence since $\D$-modules satisfy $h$-descent~\cite[Prop.\ 3.2.2]{crystals}. To compare this with the right adjoint map, note that a limit of fully faithful functors is fully faithful, so it suffices to show that $\Shv'(\bm{I} \backslash \partial\wt{\Fl}_{\bm{w}}) \subset \Shv(\bm{I} \backslash \partial\wt{\Fl}_{\bm{w}})$ and $\lim_{\bm{w}' \hra \bm{w}} \Shv'(\bm{I} \backslash \wt{\Fl}_{\bm{w}'}) \subset \lim_{\bm{w}' \hra \bm{w}} \Shv(\bm{I} \backslash \wt{\Fl}_{\bm{w}'})$ determine corresponding subsets of objects under the above equivalence. This is clear. 
	\end{rmk*}
	
	\begin{thm} \label{mainthm} 
		The functor $\mc{H}_\bullet' : \mc{P}_{I, \mr{fin}}^\triangleright \to \alg(\cat)$ from~\ref{last} is a colimit diagram. 
	\end{thm}
	\begin{proof}
		By~\cite[Lem.\ 3.2.2.6]{ha}, the forgetful functor $\alg(\cat) \to \cat$ is conservative, so it suffices to show that the induced map 
		\[
			\mr{oblv}(\colim_{\mc{P}_{I, \mr{fin}}} \mc{H}_\bullet) \to \mc{H}_I
		\]
		is an equivalence of (non-monoidal) DG-categories. In view of Corollary~\ref{interpret-cor}, it is equivalent to show that the map 
		\[
			\varphi : \colim_{\ms{Rig}^{\ms{d}}} \wt{\mc{H}}_\bullet \to \mc{H}_I
		\]
		is an equivalence of DG-categories. 
		
		Let us describe $\varphi$ in a more detailed way. First,  $\mc{H}_\bullet'$ gives a natural transformation 
		\[
			\eta : \mc{H}_\bullet \to \ul{\mc{H}}_{I, \mc{P}_{I, \mr{fin}}}
		\]
		where the target is the constant functor with value $\mc{H}_I$. Since the constructions in Proposition~\ref{preamalg-prop} and Corollary~\ref{interpret-cor} are functorial, $\eta$ induces a natural transformation 
		\[
			\wt{\eta} : \wt{\mc{H}}_\bullet \to \ul{\mc{H}}_{I, \ms{Rig}^{\ms{d}}}, 
		\]
		where the target is again a constant functor. Furthermore, replacing $\mc{H}_\bullet$ by $\ul{\mc{H}}_{I, \mc{P}_{I, \mr{fin}}}$ in the third bullet of Corollary~\ref{interpret-cor} gives 
		\[
			\colim_{\ms{Rig}^{\ms{d}}} \ul{\mc{H}}_{I, \ms{Rig}^{\ms{d}}} \simeq \colim_{\mc{P}_{I, \mr{fin}}} \ul{\mc{H}}_I \simeq \mc{H}_I. 
		\]
		Now we can identify $\varphi$ with the colimit of the natural transformation $\wt{\eta}$. 
		
		Let us interpret $\wt{\eta}$ as a lax-commutation datum in the following diagram: 
		\begin{cd}
			\ms{Rig}^{\ms{d}} \ar[r, "\wt{\mc{H}}_\bullet"] \ar[d, swap, "r"] & \cat \\
			L \ar[ru, swap, "\ul{\mc{H}}_{I, L}"] \ar[draw=none, bend left = 20]{ru}[anchor=center]{\Downarrow \wt{\eta}}
		\end{cd}
		Define a subfunctor $\mc{H}_{\le \bullet} \hra \ul{\mc{H}}_{I, L}$ as follows: for each $j = (l, \omega) \in L$, the subcategory $\mc{H}_{\le j} \hra \mc{H}_{I}$ is characterized by the requirement that the sheaf is supported on the union of Schubert varieties $\bigcup_{\omega' \in W_{\le j}} \Fl_{\omega'}$, where $W_{\le j}$ consists of all $\omega' \in W$ such that $(\ell(\omega'), \omega') \le j$ as elements of $L$. The definition of the total order on $L$ implies that $\wt{\eta}$ factors (uniquely) through this subfunctor, giving a natural transformation $\wt{\eta}'$ in the following diagram: 
		\begin{cd}
			\ms{Rig}^{\ms{d}} \ar[r, "\wt{\mc{H}}_\bullet"] \ar[d, swap, "r"] & \cat \\
			L \ar[ru, swap, "\mc{H}_{{\le}\bullet}"] \ar[draw=none, bend left = 20]{ru}[anchor=center]{\Downarrow \wt{\eta}'}
		\end{cd}
		We will show that $\wt{\eta}'$ induces an equivalence 
		\[
			\LKE_r \wt{\mc{H}}_\bullet \xra{\sim} \mc{H}_{{\le}\bullet}. 
		\]
		Then, taking the left Kan extensions of these two functors along $L \to \pt$ shows that $\varphi$ is an isomorphism, since $\colim_L \mc{H}_{{\le}\bullet} \simeq \mc{H}_I$. 
		
		By the explicit formula for left Kan extension, it suffices to prove the following: 
		\begin{itemize}
			\item For each $j \in L$, the map $\displaystyle\colim_{\ms{Rig}^{\ms{d}}_{\le j}} \wt{\mc{H}}_\bullet \to \mc{H}_{\le j}$ is an equivalence. 
		\end{itemize}
		Since $L$ is well-ordered, we may apply transfinite induction along $j \in L$.  The base case $j = (0, 1)$ is obvious. The limit case is trivial because colimits commute with colimits. For the successor case, we assume that $j$ has a predecessor and that the statement holds with $<j$ in place of $\le j$. 
		
		Since we have already verified~\ref{smooth}(3) in Corollary~\ref{blowup-cor}, we may apply the statements~\ref{smooth}(i, ii). Split into cases accordingly: 
		\begin{enumerate}[label=(\roman*)]
			\item Assume that $j$ satisfies $l = \ell(\omega)$. Then $W_{<j} = W_{\le j} \setminus \{\omega\}$. Choose any $\bm{w} \in \ms{Rig}^{\ms{d}}_j$, and consider the diagram 
			\begin{cd}
				\displaystyle\colim_{\la \ms{Emb}^{\ms{d}}_{<j} \xra{\ms{Emb}^{\ms{d}}} \bm{w} \ra} \wt{\mc{H}}_\bullet \circ \ms{tail} \ar[r] \ar[d] & \wt{\mc{H}}_{\bm{w}} \ar[d] \\
				\displaystyle\colim_{\ms{Rig}^{\ms{d}}_{<j}} \wt{\mc{H}}_\bullet \ar[r] \ar{d}[anchor=north, rotate=90]{\sim} & \displaystyle\colim_{\ms{Rig}^{\ms{d}}_{\le j}} \wt{\mc{H}}_\bullet \ar[d] \\
				\mc{H}_{<j} \ar[r] & \mc{H}_{\le j}
			\end{cd}
			where the upper square is cocartesian and comes from~\ref{smooth}(i), and the lower square comes from $\wt{\eta}'$. Lemma~\ref{blowup} and the proof of Corollary~\ref{blowup-cor} together imply that the outer square is cocartesian, so the lower square is also cocartesian. In the lower square, the left vertical map is an equivalence by the inductive hypothesis, so the right vertical map is also an equivalence, which completes the inductive step. 
			\item Assume that $j$ satisfies $l > \ell(\omega)$. Then $W_{<j} = W_{\le j}$, so $\mc{H}_{<j} \simeq \mc{H}_{\le j}$. On the other hand,~\ref{smooth}(ii) implies that $\ms{colim}_{\ms{Rig}^{\ms{d}}_{<j}} \wt{\mc{H}}_\bullet \simeq \ms{colim}_{\ms{Rig}^{\ms{d}}_{\le j}} \wt{\mc{H}}_\bullet$, so the desired statement follows directly from the inductive hypothesis. \qedhere
		\end{enumerate}
	\end{proof}
	
\newcommand\arXiv[1]{\href{http://arxiv.org/abs/#1}{\tt arXiv:#1}}

\end{document}